\documentclass[12pt,twoside]{article}

\usepackage{amsmath,amsfonts,amsthm,amssymb,amsopn,amstext,amscd}
\usepackage{enumerate}
\usepackage{graphicx}
\usepackage{float}
\usepackage{xcolor}
\usepackage{longtable}
\usepackage[font=footnotesize]{caption}
\usepackage{setspace} 
\usepackage{enumitem}
\usepackage{soul}
\usepackage{hyperref}
\allowdisplaybreaks

\newcommand{\N}{\mathbb{N}}

\newcommand{\cuad}{\begin{flushright}\vspace{-2ex}$\Box$\vspace{-2ex}\end{flushright}}

\newenvironment{Prf}[1][\unskip]{%
\par
\noindent
{\textbf{Proof of #1}}\newline
\vspace{-2ex}\noindent{}\newline}\cuad

\newenvironment{pro}[1][\unskip]{%
\par
\noindent
{\textbf{Proof.}}\newline
\vspace{-2ex}\noindent{}\newline}\cuad

\newenvironment{abs}
{
\begin{center} {\Large\textbf{Abstract}}
\end{center}\begin{quote}\small}
{\end{quote}
}

\newtheorem{thm}{Theorem}
\newtheorem{prop}{Proposition}

\newtheorem{cor}{Corollary}
\newtheorem{lem}{Lemma}
\newtheorem{remark}{Remark}

\newtheorem{defi}{Definition}

\newcommand{\vc}[1]{\boldsymbol{#1}}
\usepackage{dsfont}
\newcommand{\ind}[1]{\mathds{1}_{\{#1\}}}

\newcommand{\bm}[1]{\mbox{\boldmath $#1$}}
\newcommand{\nc}{\newcommand}
\nc{\mbN}{\mathbb N}
\nc{\mbP}{\mathbb P}
\nc{\mbR}{\mathbb R}
\nc{\mbE}{\mathbb E}
\nc{\mbV}{\mathsf{Var}}

\begin{document}
\pagenumbering{arabic}

%
%
%
%
%
%
%
%
%
%
%
%
%

\newpage

\pagenumbering{arabic}

\title{Parameter estimation in
branching processes with almost sure extinction}
\author{P. Braunsteins\footnote{School of Mathematics and Statistics, University of Melbourne, email: \url{p.braunsteins@unimelb.edu.au}} \footnote{Korteweg-de Vries Instituut,  University of Amsterdam, email: \url{p.t.m.braunsteins@uva.nl}, ORCID: 0000-0003-1864-0703}, S. Hautphenne\footnote{School of Mathematics and Statistics, University of Melbourne, email: \url{sophiemh@unimelb.edu.au}, ORCID: 0000-0002-8361-1901}, C. Minuesa\footnote{{Department of Mathematics, Autonomous University of Madrid, email: \url{carmen.minuesa@uam.es}, ORCID: 0000-0002-8858-3145.
Note: all authors contributed equally to this work.}}}
\maketitle
\date{ }

\begin{abs} 
We consider population-size-dependent branching processes (PSDBPs) which eventually become extinct with probability one. 
For these processes, we derive maximum likelihood estimators for the mean number of offspring born to individuals when the current population size is $z\geq 1$.
As is standard in branching process theory, an asymptotic analysis of the estimators requires us to condition on non-extinction up to a finite generation $n$ and let $n\to\infty$; however, because the processes become extinct with probability one, we are able to demonstrate that our estimators do not satisfy the classical consistency property ($C$-consistency). 
This leads us to define the concept of \emph{$Q$-consistency}, and we prove that our estimators are $Q$-consistent and asymptotically normal.
To investigate the circumstances in which a $C$-consistent estimator is preferable to a $Q$-consistent estimator, we then provide two $C$-consistent estimators for subcritical Galton--Watson branching processes.
Our results rely on a combination of linear operator theory, coupling arguments, and martingale methods.

\end{abs}

\section{Introduction}

Branching processes are the primary tool used to model populations that evolve randomly in time.
They have been used successfully to illuminate real-world problems arising in many areas, in particular in biology and conservation ecology \cite{haccou,jager,kimmel2002branching}, which are our primary focus here.
For these reasons, parameter estimation in branching processes has received significant attention; we refer the reader to the survey by Yanev \cite{yanev08}, and to Guttorp \cite{Guttorp} for a book-length treatment.

The simplest discrete-time branching process is the Galton--Watson (GW) process. In this process, at every generation, individuals reproduce independently  according to a common offspring distribution $\xi$ with mean $m:=\mbE(\xi)$. The population size at generation $n$, denoted by $Z_n$, satisfies the well-known recursion $Z_n=\sum_{i=1}^{Z_{n-1}}\xi_{n,i}$ for $n\geq1$, where $\{ \xi_{n,i} \}_{ n,i \in \mathbb{N}_0}$ are independent copies of $\xi$.
Several estimators have been developed to estimate the mean offspring $m$ based on the observation of the successive population sizes $Z_0,Z_1, \ldots,Z_n$. The most efficient estimator for $m$ is the maximum likelihood estimator (MLE) introduced by Harris in \cite{harris48},
\begin{equation}\label{mlem}\hat{m}_n:= \dfrac{\sum_{i=1}^{n} Z_i}{\sum_{i=1}^{n}Z_{i-1}}, \end{equation}which can be identified as the ratio between the total cumulative number of children over the total cumulative number of parents in the observed sample.
Harris showed that, in the supercritical case $m>1$, on the set of non-extinction $\{Z_n\rightarrow\infty\}$, this estimator is consistent. In other words, for any initial state $i\geq 1$ and any $\varepsilon>0$, 
\begin{eqnarray}\label{Ccons} \lim_{n\rightarrow\infty} \mbP_i(|\hat{m}_n-{m}|>\varepsilon\,|\,Z_n>0)=0,\end{eqnarray}where $\mbP_i(\cdot)$ stands for $\mbP(\cdot | Z_0=i)$; this property is also called \emph{$C$-consistency} of the estimator \cite{Pakes-1975}. Conditional asymptotic normality of $\hat{m}_n$ in the supercritical case was proved later by Dion~\cite{dion74}.


Because $\mbE(Z_n\,|\,Z_0)=Z_0\, m^n$, GW processes grow exponentially.  Endless exponential growth is however rarely observed in real biological populations, which instead tend to exhibit \emph{logistic growth}. Indeed, due to the presence of competition between individuals, as the population size increases, its growth starts to slow down until the population reaches a threshold called the \emph{carrying capacity}, which corresponds to the maximum population size that the habitat can support. The population then fluctuates around this carrying capacity for a long period of time, before eventually becoming extinct \cite{hamza2016establishment,jk11,jz20}. 
Of all species that have ever existed it is estimated that over 99.9\% are now extinct \cite{raup}, so it does indeed appear that most populations face eventual extinction; moreover, under standard conditions, the background extinction rate is of the order of one in millions of years \cite{newman97}, so it also appears that these populations tend to survive for a long period of time ---often much longer than typical observation/census periods.

Population-size-dependent branching processes (PSDBPs) are suitable models for populations that display logistic growth. They are defined similarly to GW processes, the only difference being that, at each generation, the offspring distribution $\xi(z)$ now depends on the current population size $z$ ($z\geq 1$). Accordingly, the mean offspring $m(z):=\mbE[\xi(z)]$ also becomes a function of the current population size $z$. The carrying capacity is then defined as the threshold value $K$ such that $m(z)>1$ as long as $z< K$, and
$m(z)<1$ when $z>K$. In this paper we estimate $m(z)$ for any population size $z$, based on the observation of successive generation sizes. A natural analogue of \eqref{mlem}  is the estimator
\begin{equation}\label{mlemz}
\hat{m}_{n}(z):=\frac{\sum_{i=1}^{n} Z_{i}\,\ind{Z_{i-1}=z}}{z\sum_{i=1}^{n}\ind{Z_{i-1}=z}}.
\end{equation}In Proposition \ref{prop:mle-m}, we prove  that $\hat{m}_{n}(z)$ is the MLE for $m(z)$. In order to establish consistency and asymptotic normality of this estimator, we need \textit{infinitely long} trajectories of the process $\{Z_n\}$, and the standard way to proceed is to condition on survival of the process. However, we are in a unique setting where trajectories are typically very long but have no chance to survive forever. Conditioning on survival then pushes the sample paths away from the absorbing  state 0 (extinction), and therefore induces a bias in the estimates. As a consequence, conditional on $Z_n>0$, the estimator $\hat{m}_{n}(z)$ does  \textit{not} converge to the true mean offspring $m(z)$ but to a different limit that we denote by $m^{\uparrow}(z)$. We show that $m^{\uparrow}(z)$ has a probabilistic interpretation in terms of the \emph{$Q$-process} $\{Z^{\uparrow}_n\}_{n\geq0}$, or Doob $h$-transform of the original PSDBP $\{Z_n\}_{n\geq0}$.  The $Q$-process is a non-absorbing Markov chain whose behaviour is very similar to the original PSDBP except in the vicinity of the absorbing state, and which can be interpreted as the original process conditioned on not becoming extinct in the distant future. The limit $m^{\uparrow}(z)$ of the estimator $\hat{m}_{n}(z)$ is then \emph{the equivalent} of the mean offspring  $m(z)$ in the $Q$-process. More formally, in Theorem \ref{thm:psdbp-Qcons and Qnorm} we show that for any initial state $i\geq 1$,  any $z\geq 1$, and any $\varepsilon>0$, 
\begin{eqnarray}\label{Qcons} \lim_{n\rightarrow\infty} \mbP_i(|\hat{m}_n(z)-{m}^{\uparrow}(z)|>\varepsilon\,|\,Z_n>0)=0.\end{eqnarray}We name this property \emph{$Q$-consistency} of the estimator $\hat{m}_n(z)$. We complement this result by showing conditional asymptotic normality of the estimator. 

In summary, there are now two concepts of consistency: $C$-consistency, when, conditional on survival, the estimator converges to the true parameter, and $Q$-consistency, when, conditional on survival, the estimator converges to the equivalent of the true parameter in the $Q$-process. Since the trajectories of PSDBPs that do not become extinct quickly generally survive for a very long time, the functions $m(z)$ and $m^{\uparrow}(z)$ are  close to each other in the vicinity of the carrying capacity, especially if the latter is large. Thus, in practice,  the $Q$-consistent estimator $\hat{m}_n(z)$ (which is the MLE) is often satisfactory. However, a natural question follows: Are there situations where a $C$-consistent estimator is preferable to a $Q$-consistent one? We often study endangered populations \textit{because} they are still alive. For these populations, we therefore have to think about the observed population sizes as being generated under the condition $Z_n>0$, which biases the sample. And in this context, where important conservation measures may need to be taken, we may favour an estimator which converges to the true parameter value.

To investigate this question, we first need to derive a $C$-consistent estimator. To do so, we return to the GW setting and focus on the subcritical regime, where extinction occurs rapidly. 
Scant attention has been paid in the literature to statistical inference of subcritical GW processes when the number of ancestors is fixed and in absence of immigration \cite{yanev08}. 
We show that, conditional on $Z_n>0$, and regardless the value of $m<1$, the MLE $\hat{m}_n$ always converges to one (Theorem~\ref{thm:GW-Qcons and Qnorm}); this result can be interpreted as $Q$-consistency of  $\hat{m}_n$ (Corollary~\ref{cor:GW-Qcons}). 
We then use the probabilistic interpretation of the $Q$-process associated with a GW process to derive two $C$-consistent (and asymptotically normal) estimators for the mean offspring (Propositions \ref{Ccons1} and \ref{Ccons2}). The first one relies on a stronger assumption on the offspring distribution than the second one but is more efficient. To the best of our knowledge, these are the first $C$-consistent estimators for subcritical branching processes, and we point out that they are also $C$-consistent in the critical and supercritical regimes. We compare the $C$-consistent estimators with the ($Q$-consistent) MLE in a numerical example.

We prove Theorem \ref{thm:psdbp-Qcons and Qnorm} (for PSDBPs) and Propositions \ref{Ccons1} and \ref{Ccons2} (for GW processes) using a unified approach. 
More specifically, we form a \textit{MEXIT} coupling (or \textit{maximal exit time} coupling \cite{ernst2019mexit}) of the original branching process \textit{conditional on $Z_n>0$} (which results in a non-homogeneous Markov chain) and its corresponding $Q$-process (which is a more tractable homogeneous Markov chain). Roughly speaking, in the MEXIT coupling, the probability that the trajectories of the two processes stick together until any generation $n-k$ is maximised.
We use the linear operator theory approach of Gosselin~\cite{Gosselin-2001} to 
bound the probability that the processes uncouple by generation $n-k$ (Theorem \ref{coupling1}).
The idea  is then to consider the two processes to be the same at least up to generation $n-k$ (where $k\sim \log n$), and to show that ``nothing bad'' can happen for the remaining $k$ generations. Studying the asymptotic properties of the estimators conditional on $Z_n>0$ now reduces to studying these properties in the $Q$-process, and we complete our argument by applying the martingale central limit theorem. The bounds in Theorem \ref{coupling1} also provide us with an efficient method of approximate simulation of long non-extinct trajectories 
with a controlled error. 

We highlight that, while this paper focuses on parameter estimation in branching processes with almost sure extinction, the ideas extend well beyond this scope, namely to other absorbing Markov chains for which parameter estimation is of interest. Examples include estimation of the transmission rate/probability in stochastic epidemic models where the number of infected individuals eventually reaches zero with probability one (such as \textit{SIS} models). In addition, when we consider PSDBPs, we focus on the estimators $\hat{m}_n(z)$ for the mean offspring at population sizes $z\geq 1$. The asymptotic properties of these $Q$-consistent estimators (Theorem \ref{thm:psdbp-Qcons and Qnorm}) have further implications: in a subsequent paper, we use these properties to derive $C$-consistent estimators for parameters in PSDBPs whose mean offspring has a parametric form, such as the Beverton-Holt and the Ricker models (defined in \eqref{bh} and \eqref{ricker}). 


The paper is organised as follows. In Section 2 we give some preliminary results:  we provide  background on PSDBPs (Section 2.1), we introduce the MLE $\hat{m}_{n}(z)$ for $m(z)$ (Section 2.2), and we define the $Q$-process $\{Z^{\uparrow}_n\}$ associated with a branching process $\{Z_n\}$ (Section 2.3). In Section 3, we present our main results: we provide the asymptotic properties of $\hat{m}_{n}(z)$ and define the concept of $Q$-consistency (Section 3.1), we present $C$-consistent estimators for subcritical GW processes (Section 3.2), and we provide properties of the MEXIT coupling (Section 3.3).    
The proofs of our results are gathered in Section \ref{proofs}.

Throughout the paper we use the following notation: we let $\N_0:=\N\cup\{0\}$, we let $\vc 1$  be a column vector with all components equal to 1, we let $\vc e_i$ be a column vector whose $i$-th element is equal to 1 and the remaining elements are 0, and we let $\mathds{1}_A$ be the indicator function of the set $A$. In addition, we sometimes use the shorthand notation $X\,|\,B$ to  denote a random variable $X$ conditional on the event $B$.


\section{Preliminaries}\label{sec:model}

\subsection{Population-size-dependent branching processes}\label{sec:psdbp}
A discrete-time \emph{population-size-dependent branching process} (PSDBP) is a 
 process in which, at each generation, individuals reproduce independently and give birth to  a random number of offspring $\xi(z)$ that depends on the population size $z\geq 1$ in that generation. We denote the offspring distribution at population size $z$ by $\boldsymbol{p}(z)=(p_k(z))_{k\in\N_0}$, where $p_k(z):=\mbP(\xi(z)=k)$, and we assume that $p_0(z)>0$ and $p_0(z)+p_1(z)<1$ for each $z$. 
We define the offspring mean and variance functions as $m(z):=\mbE[\xi(z)]$ and $\sigma^2(z):=\mbV[\xi(z)]$, respectively, and assume that both  are finite for each $z$.  

Let $Z_n$ represent the population size at generation $n$;
the process  $\{Z_n\}_{n\in\N_0}$ is characterised by the recursive equation
%
%
%
\begin{equation}\label{def:process}
Z_0=N,\quad Z_{n+1}=\sum_{i=1}^{Z_n} \xi_{n,i}(Z_n),\quad n\in\N_0,
\end{equation}
where $N$ is a positive integer, and where for each $z\geq 1$, $\{\xi_{n,i}(z):n\in\N_0, i=1,\ldots,z\}$ is a family of independent random variables with the same distribution as $\xi(z)$. The empty sum in \eqref{def:process} is taken to be 0.  The process $\{Z_n\}$ is therefore a time-homogeneous Markov chain with the absorbing state 0.
The conditional moments of the process are given by:
\begin{equation}\label{eq:prop-moment-proc}
\mbE(Z_{n+1}|Z_n)=Z_n m(Z_n),\quad \text{ and }\quad\mbV(Z_{n+1}|Z_n)=Z_n \sigma^2(Z_n).
\end{equation}


If the offspring distribution is independent of the current population size $z$, then $\{Z_n\}$ reduces to the standard Galton--Watson (GW) process with offspring distribution $\boldsymbol{p}$, offspring mean $m$, and offspring variance $\sigma^2$. 

\medskip

Two well known models for PSDBPs with a \textit{carrying capacity} are the Beverton-Holt model, in which 
\begin{equation}\label{bh}m(z)=\frac{2K}{K+z},\qquad K>0,\end{equation} and the Ricker model, in which 
\begin{equation}\label{ricker}m(z)=r^{(1-z/K)},\qquad r>1, K>0;\end{equation}see for instance \cite{hognas}. In these models, $K$ represents the carrying capacity: if $z<K$ then $m(z)>1$, whereas if $z>K$ then $m(z)<1$.
In the Beverton-Holt model, the expected next generation size $\mbE(Z_{n+1}|Z_n=z)=z m(z)$ is increasing as a function of $z$, while in the Ricker model, it approaches zero as $z$ becomes large (significant over-population may reduce the size of the next generation dramatically).

\begin{figure}[t]
\centering\includegraphics[width=1\textwidth]{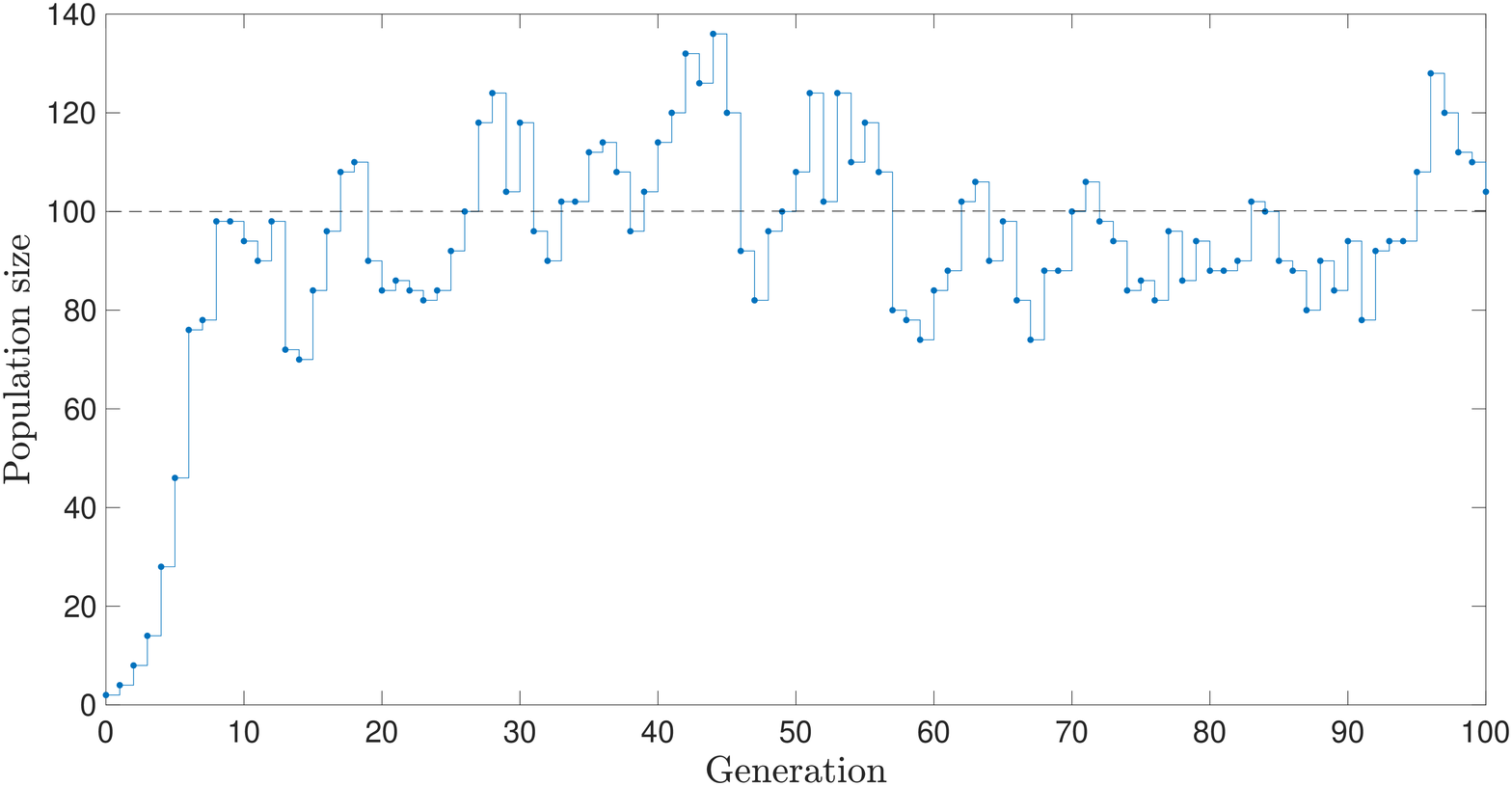}
\caption{\label{bhf1}Trajectory of the Beverton-Holt model  with carrying capacity $K=100$, $p_2(z)=m(z)/2$, and $p_0(z)=1-p_2(z)$ (binary splitting).}
\end{figure}

PSDBPs with a fixed carrying capacity eventually become extinct with probability one, that is, $\mbP_i(Z_n\rightarrow0)=1$ for any initial population size $i$; see for instance \cite{jager92}. If the population enters the vicinity of the carrying capacity, it typically lingers around that threshold value for a very long time \cite{hamza2016establishment,jk11}; see Figure \ref{bhf1} for an illustration.

\subsection{Maximum likelihood estimation of $m(z)$}\label{sec:mle-pk}

In this section we study the maximum likelihood estimator (MLE) for the mean offspring at population size $z$, $m(z)$, based on the sample $\mathcal{Z}_n:=\{Z_0,\ldots,Z_{n}\}$ which contains the successive population sizes up to generation $n$. Except in the binary splitting case,  MLEs for the full distribution $\vc p(z)$ and for the variance $\sigma^2(z)$ require a more refined sample, as discussed in Section \ref{ape:mle-pk}.

For $z,n\in\N$, we define
$$j_n(z):=\sum_{i=1}^{n}\ind{Z_{i-1}=z},$$
which is the number of times the population size equals $z$ during the first $n$ generations.
The next proposition states that the natural estimator for $m(z)$ given in \eqref{mlemz} is the MLE. This result is part of Corollary \ref{cor:mle-m-sigma}, which is stated and proved in Appendix \ref{ape:mle-pk}.
\begin{prop}\label{prop:mle-m}For every $z\in\N$, the MLE of $m(z)$ on the set $\{j_n(z)>0\}$ based on the sample $\mathcal{Z}_n$ is
\begin{equation*}\label{eq:mle-m}
\begin{aligned}
\hat{m}_{n}(z)&:=\frac{\sum_{i=1}^{n} Z_{i}\ind{Z_{i-1}=z}}{z\sum_{i=1}^{n}\ind{Z_{i-1}=z}}=\frac{\sum_{i=1}^{n} Z_{i}\ind{Z_{i-1}=z}}{z\,j_n(z)}.
\end{aligned}
\end{equation*}
\end{prop}

In the next proposition we derive the conditional mean of $\hat{m}_{n}(z)$ given \mbox{$j_n(z)>0$}.

%
\begin{prop}\label{prop:moments-mle}
For each  $z\in\N$ and $n\in\N$,
\begin{equation}\mbE\left(\hat{m}_{n}(z)\,|\,j_n(z)>0\right)=m(z)\mbE\left(\frac{\ind{Z_{n-1}=z}}{j_n(z)}\,\Big|\, j_n(z)>0\right)+\mbE\left(\frac{\sum_{i=1}^{n-1}Z_{i}\ind{Z_{i-1}=z}}{zj_n(z)}\,\Big|\, j_n(z)>0\right).\label{prop:mean-mle-m}\end{equation}
\end{prop}
Proposition \ref{prop:moments-mle} is the PSDBP counterpart of Proposition 2.4 in \cite{Guttorp} (for GW processes) and Proposition~4.2 in \cite{art-2004b}  (for controlled branching processes). 
\begin{remark}\label{rem:biasMLE}
If the variables $Z_1,\ldots, Z_{n-1}$ were independent of $j_n(z)$, then the estimator $\hat{m}_{n}(z)$ would be unbiased. However, these variables are clearly not independent of $j_n(z)$: indeed, if $j_n(z)=n$, then $Z_j=z$ for all $j=1,\ldots, n-1$. The variable $Z_n$ is, however, conditionally independent of $j_n(z)$, given $Z_{n-1}=z$.  

\end{remark}

\begin{figure}[t]
\centering\includegraphics[width=1\textwidth]{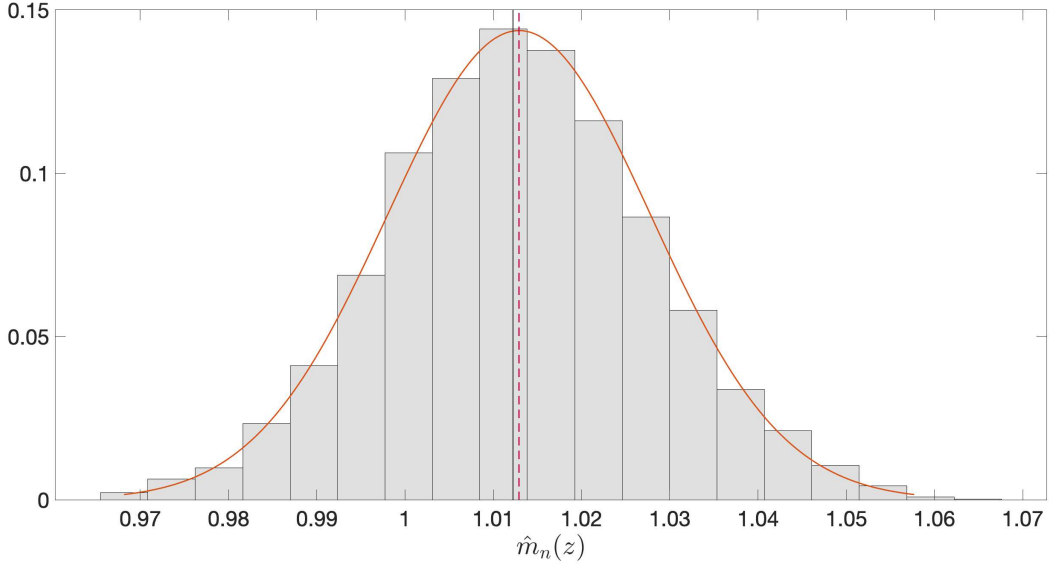}
\caption{\label{rick3a}Binary splitting Ricker model \eqref{ricker} with $r=1.2$ and carrying capacity $K=30$. Histogram of the estimates $\hat{m}_n(z)$ for $n=2000$ and \textbf{$z=28$} ($n u_zv_z=159.94$), based on $5000$ simulated non-extinct trajectories. Black vertical line: true $m(z)$ (1.0122); Red dashed line: true $m^{\uparrow}(z)$ (1.0129); Blue dashed line: empirical mean of $\hat{m}_n(z)$ (1.0129) (note that the red and blue lines are superimposed); Red curve: theoretical asymptotic normal distribution given in Theorem \ref{thm:psdbp-Qcons and Qnorm}.}
\end{figure}

To investigate the behaviour of the MLE $\hat{m}_{n}(z)$ for large values of $n$, we first look at an example. We consider the Ricker model \eqref{ricker} with $r=1.2$ and $K=30$, and with \emph{binary splitting}, in which individuals have two offspring with probability $p_2(z)=m(z)/2$ and no offspring with probability $p_0(z)=1-p_2(z)$.
Figure \ref{rick3a} shows a histogram of the estimates $\hat{m}_{n}(z)$ for $z=28$, based on 5000 simulated non-extinct trajectories of length $n=2000$. The true value $m(28)=1.0122$ is represented by the black vertical line, and the empirical mean of the estimates, equal to 1.0129, is represented by the dashed line. The graph indicates that, conditionally on $Z_n>0$, $\hat{m}_{n}(z)$ has an asymptotic normal distribution.

To prove asymptotic normality, $n$ needs to be taken to infinity, which requires infinitely long non-extinct trajectories; this means that we need to condition on survival of the process. However, we are in a setting where the process becomes extinct with probability one. Conditioning on survival in this case pushes the sample paths away from extinction, which results in overestimating $m(z)$, especially for low values of $z$. As a consequence, the MLE does \emph{not} converge to the true mean offspring, but to a different value, which we denote by $m^{\uparrow}(z)$. In our example, since $z=28$ is close to the carrying capacity, the probability of extinction in the short term from state $z$ is relatively low, so conditioning on survival does not significantly `modify' the transition probabilities from $z$ in the original process, hence $m(28)\approx m^{\uparrow}(28)=1.0129$. This is formalised in Sections \ref{subsec:Q-process} and \ref{sec:ass_prop}, in which we provide a probabilistic interpretation of $m^{\uparrow}(z)$.

\subsection{Branching processes conditioned on non-extinction}\label{subsec:Q-process}


In this section we describe the effect of the condition $Z_n>0$ on the transition probabilities of $\{Z_n\}$.
We denote by $Q$ the sub-stochastic transition probability matrix of $\{Z_n\}$ restricted to the transient states $\{1,2,\ldots\}$,  and we make the following regularity assumptions:

\begin{enumerate}[label=(A\arabic*),ref=(A\arabic*),start=1]
\item There exists $z\in\N$ and $n\geq 1$ such that $(Q^n)_{zz}>0$.\label{cond:irreducible-PSDBP}
\item $\limsup_{z\to\infty} m(z)<1$.\label{cond:lim-sup-m-z}
\item For each $\nu\in\N$, $\sup_{z\in\N} \mbE[\xi(z)^\nu]<\infty$\label{cond:bounded-moments-xi}.
\end{enumerate}
Under Assumptions \ref{cond:lim-sup-m-z} and \ref{cond:bounded-moments-xi}, the process becomes extinct with probability one for any initial distribution $\vc \pi$ (see \cite[Proposition 3.1]{Gosselin-2001}), that is, $\mbP_\pi(Z_n\to 0)=1$, where $\mbP_\pi(\cdot)$ denotes the probability measure given the initial distribution $\vc \pi$. Moreover, following \cite[Theorem 4.1]{Gosselin-2001}, under Assumptions \ref{cond:irreducible-PSDBP}--\ref{cond:bounded-moments-xi}, for any $n\geq1$ we can write $Q^n$ uniquely as
\begin{equation}\label{PF1}
Q^n=\rho^n \vc v\vc u^\top +S^n, 
\end{equation} 
where $\rho:=\rho(Q)=\lim_{n\rightarrow\infty} (Q^n)_{ij}^{1/n}$ is the convergence norm of $Q$, $\vc u$ and $\vc v$ are real strictly positive column vectors such that
\begin{eqnarray}\label{uv}\vc u^\top Q=\rho \vc u^\top, \quad Q\vc v=\rho \vc v, \quad \vc u^\top\vc 1=1,\quad \textrm{and}\quad \vc u^\top\vc v=1,\end{eqnarray}
and $S$ is a matrix whose convergence norm is such that $\rho(S)<\rho$, so that
%
\begin{equation}\label{Qn}
Q^n\sim\rho^n \vc v\vc u^\top, \quad\textrm{as $n\rightarrow\infty$};
\end{equation}see also Lemma \ref{Goss2} in Section \ref{lot}.
The vector $\vc u$ corresponds to the \emph{quasi-stationary} (or \emph{quasi-limiting}) distribution of $\{Z_n\}$, and the vector $\vc v$ records the relative ``strength'' of each state. Indeed, using \eqref{Qn},
$$\lim_{n\rightarrow\infty}\mbP_i(Z_n=j\,|\,Z_{n}>0)=\lim_{n\rightarrow\infty}\dfrac{\vc e_i^\top Q^n \vc e_j}{\vc e_i^\top Q^{n}\vc 1}=\lim_{n\to\infty}\dfrac{\rho^n v_i u_j}{\rho^n v_i}=u_j,\qquad j\geq 1,$$and
$$\lim_{n\rightarrow\infty}\dfrac{\mbP_j(Z_{n}>0)}{\mbP_i(Z_{n}>0)}=\lim_{n\rightarrow\infty}\dfrac{\vc e_j^\top Q^n \vc 1}{\vc e_i^\top Q^{n}\vc 1}=\dfrac{v_j}{v_i},\qquad i,j\geq 1,$$where recall that $\mbP_i(\cdot)$ stands for $\mbP(\cdot | Z_0=i)$. 
In addition, by  \cite[Theorem 3.1]{Gosselin-2001}, the vectors $\vc u$ and $\vc v$ satisfy
\begin{eqnarray}\label{uv2}
\sum_{i\geq 1} i \,u_i<\infty, \quad \textrm{and}\quad \sup_{i\in\N}\frac{v_i}{i}<\infty.
\end{eqnarray}


 For every $n\geq0$ fixed, the process $\{Z_\ell\}_{0\leq \ell\leq n}$ conditioned on $Z_{n}>0$ is a \textit{time-inhomogeneous} Markov chain that we denote by $\{Z_\ell^{(n)}\}_{0\leq \ell\leq n}$ and whose one-step transition probabilities are
\begin{eqnarray}\nonumber
\mbP(Z_{\ell+1}^{(n)}=j\,|\,Z_{\ell}^{(n)}=i)&:=&\mbP(Z_{\ell+1}=j\,|\,Z_{\ell}=i,\;Z_{n}>0)\\\label{condProb}&=&Q_{ij}\,\dfrac{\vc e_j^\top Q^{n-\ell-1} \vc 1}{\vc e_i^\top Q^{n-\ell}\vc 1}\quad \text{ for }i,j\geq 1.
\end{eqnarray}
Taking the limit as $n\rightarrow\infty$ in \eqref{condProb} and using \eqref{Qn} leads to homogeneous transition probabilities: 
\begin{eqnarray*}
\mbP(Z_{\ell+1}^\uparrow=j\,|\,Z_{\ell}^\uparrow=i) &:=& \lim_{n\rightarrow\infty}\,\mbP(Z_{\ell+1}^{(n)}=j\,|\,Z_{\ell}^{(n)}=i) \\ &=& \lim_{n\rightarrow\infty}Q_{ij}\,\dfrac{\vc e_j^\top \rho^{n-\ell-1} \,\vc v }{\vc e_i^\top \rho^{n-\ell} \,\vc v } \\
&=&Q_{ij}\,\dfrac{v_j}{\rho v_i}.
\end{eqnarray*}
These transition probabilities define a new positive-recurrent \emph{time-homogeneous} Markov chain $\{Z^{\uparrow}_\ell\}_{ \ell\geq 0}$, which we refer to as the \emph{$Q$-process} associated with $\{Z_n\}$ (following the terminology in \cite{athreya}). The $Q$-process $\{Z^{\uparrow}_\ell\}$ can be interpreted as the original process $\{Z_n\}$ conditioned on not being extinct in the distant future.
 The $n$-step transition probabilities of the $Q$-process are given by
$$({Q^{\uparrow}}^n)_{ij}:=\mbP(Z^{\uparrow}_n=j\,|\,Z^{\uparrow}_0=i)=(Q^n)_{ij} \dfrac{v_j}{\rho^n v_i },\quad \text{ for }n,i,j\geq 1.$$
The stationary distribution of the $Q$-process follows from \eqref{Qn}:
\begin{eqnarray*}
\lim_{n\rightarrow\infty} \mbP(Z^{\uparrow}_n=j\,|\,Z^{\uparrow}_0=i) &=& u_j v_j,\quad \text{ for }i, j\geq 1. 
\end{eqnarray*}




\section{Main results}
\subsection{Asymptotic properties of the MLE $\hat{m}_{n}(z)$}\label{sec:ass_prop}

Now that we have introduced the $Q$-process, we are in a position to give a probabilistic interpretation to the limit  $m^\uparrow(z)$ that we observed in Section \ref{sec:mle-pk}. The mean offspring at population size $z$ can be written in terms of the transition probabilities of $\{Z_n\}$ as the normalised mean next state from $z$:
\begin{equation}\label{fmz}m(z)=z^{-1} \sum_{j\geq 1} j\, Q_{zj}.\end{equation} The function $m^\uparrow(z)$ is the quantity equivalent to $m(z)$ in the $Q$-process, that is,
\begin{equation}\label{fmhz}{m}^{\uparrow}(z):=z^{-1} \sum_{j\geq 1} j\, Q^{\uparrow}_{zj},\qquad \text{where}\quad Q^{\uparrow}_{zj}=Q_{zj} \dfrac{v_j}{\rho v_z }.\end{equation}In the next theorem, we formally prove that, conditional on $Z_n>0$, the estimator $\hat{m}_{n}(z)$ converges to ${m}^{\uparrow}(z)$, that is, to the analogue of the mean offspring $m(z)$ in the $Q$-process; we also show conditional asymptotic normality.
 Similar to the normalised mean ${m}^{\uparrow}(z)$, we define the normalised variance of the next state from $z$ in the $Q$-process as
\begin{equation}\label{sig2b}{\sigma^2}^{\uparrow}(z)=\dfrac{\sum_{k=1}^\infty k^2 Q^{\uparrow}_{zk}}{z^2}-({m}^{\uparrow}(z))^2.\end{equation}
\begin{thm}\label{thm:psdbp-Qcons and Qnorm}Under Assumptions \ref{cond:irreducible-PSDBP}--\ref{cond:bounded-moments-xi}, for any initial state $i\geq 1$ and every $z\geq 1$, the MLE  $\hat{m}_n(z)$ for $m(z)$ in a PSDBP  satisfies,  for any $\varepsilon>0$, 
\begin{equation}\label{mz_Qcons}\lim_{n\rightarrow\infty} \mbP_i(|\hat{m}_n(z)-{m}^{\uparrow}(z)|>\varepsilon\,|\,Z_n>0)=0,\end{equation} and for any $x\in \mathbb{R}$,
\begin{equation}\label{mz_Qnorm}\lim_{n\rightarrow\infty}\mbP_i(\{n \,u_z v_z/{\sigma^2}^{\uparrow}(z)\}^{1/2} \left(\hat{m}_n(z)-{m}^{\uparrow}(z)\right)\leq x\,|\,Z_n>0)=\Phi(x),\end{equation}
where $\Phi(x)$ is the distribution function of a standard normal random variable, and where ${m}^{\uparrow}(z)$ and ${\sigma^2}^{\uparrow}(z)$ are finite \mbox{for any $z$}. 
In addition, for any pair $(z_1,z_2)\in\N^2$ with $z_1\neq z_2$, the two normalised estimators $$\sqrt{n}\left(\hat{m}_n(z_1)-{m}^{\uparrow}(z_1)\right)\quad \text{and}\quad\sqrt{n} \left(\hat{m}_n(z_2)-{m}^{\uparrow}(z_2)\right)$$ are asymptotically uncorrelated. 
\end{thm}
Theorem \ref{thm:psdbp-Qcons and Qnorm} leads us to define \emph{Q-consistency} of an estimator in an absorbing Markov chain $\{Z_n\}$ whose $Q$-process $\{Z_n^\uparrow\}$ is positive recurrent.
For $n\geq 0$, let $\hat{\theta}_n:=S_n(Z_0,\ldots, Z_n)$ be an estimator for a quantity $\theta$ in  $\{Z_n\}$, where $S_n(Z_0,\ldots, Z_n)$ is a measurable function of the first $n+1$ states of the process. For example, $\hat{\theta}_n=\hat{m}_n(z)$ for a given $z\geq1$ in a PSDBP.
\begin{defi}[$Q$-consistency]\label{def:Q-}The estimator $\hat{\theta}_n$ is  \emph{$Q$-consistent} for $\theta$ if it satisfies the following two conditions:
\begin{itemize}
\item[\emph{(i)}] If $\hat{\theta}_n^{\uparrow}:=S_n(Z_0^{\uparrow},\ldots, Z_n^{\uparrow})$, then the random sequence $\{\hat{\theta}_n^{\uparrow}\}_{n\geq0}$ converges to a constant $\theta^{\uparrow}=f(Q^\uparrow)$ with probability one, where $f(\cdot)$ is a non-constant function of the transition probabilities $Q^\uparrow_{ij}$ of $\{Z^\uparrow_n\}$ such that $f(Q)=\theta$.
\item[\emph{(ii)}] \emph{Conditional on $Z_n>0$}, for any $i\geq1$, the random sequence $\{\hat{\theta}_n\}_{n\geq0}$ converges in probability to $\theta^{\uparrow}$, that is, for any $\varepsilon>0$,
$$\lim_{n\rightarrow\infty} \mbP_i(|\hat{\theta}_n-\theta^{\uparrow}|>\varepsilon\,|\,Z_n>0)=0.$$
\end{itemize}
\end{defi}

By letting the function $f(\cdot)$ in Definition \ref{def:Q-} (i) take the specific form suggested by \eqref{fmz} and \eqref{fmhz}, we obtain the following corollary of Theorem~\ref{thm:psdbp-Qcons and Qnorm}.

\begin{cor}\label{cor:psdbp-Qcons and Qnorm}Under the assumptions of Theorem \ref{thm:psdbp-Qcons and Qnorm},  for every $z\geq1$, $\hat{m}_{n}(z)$ is a Q-consistent estimator for $m(z)$.
\end{cor}

\begin{figure}[t]
\centering\includegraphics[width=1\textwidth]{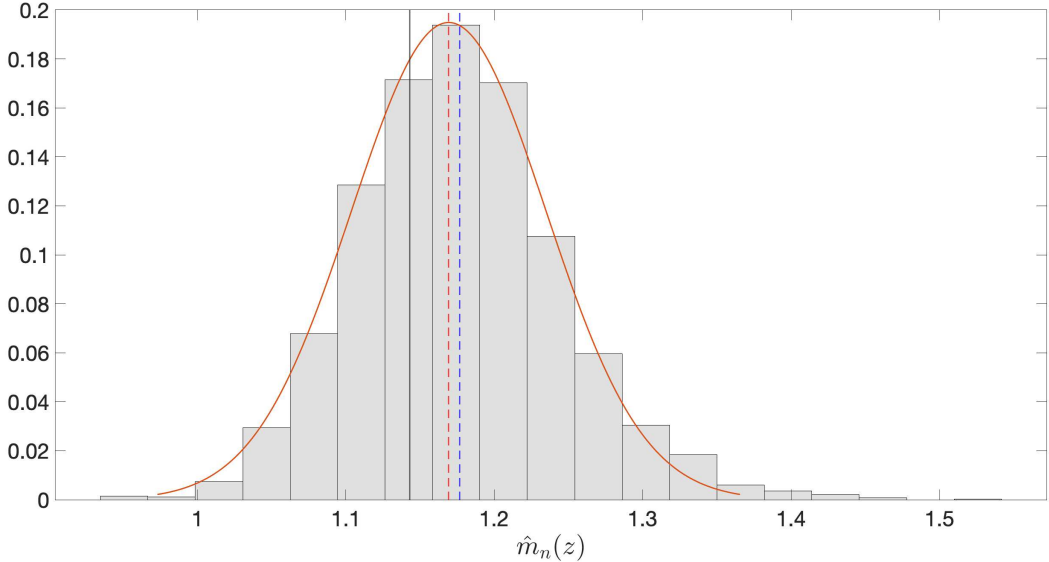}
\caption{\label{rick3b}Binary splitting Ricker model \eqref{ricker} with $r=1.2$ and carrying capacity $K=30$. Histogram of the estimates $\hat{m}_n(z)$ for $n=2000$ and \textbf{$z=8$} ($n u_zv_z=26.63$), based on $5000$ simulated non-extinct trajectories. Black vertical line: true $m(z)$ (1.1431); Red dashed line: true $m^{\uparrow}(z)$ (1.1693); Blue dashed line: empirical mean of $\hat{m}_n(z)$ (1.1768); Red curve: theoretical asymptotic normal distribution.}
\end{figure}

\begin{figure}[t]
\centering\includegraphics[width=1\textwidth]{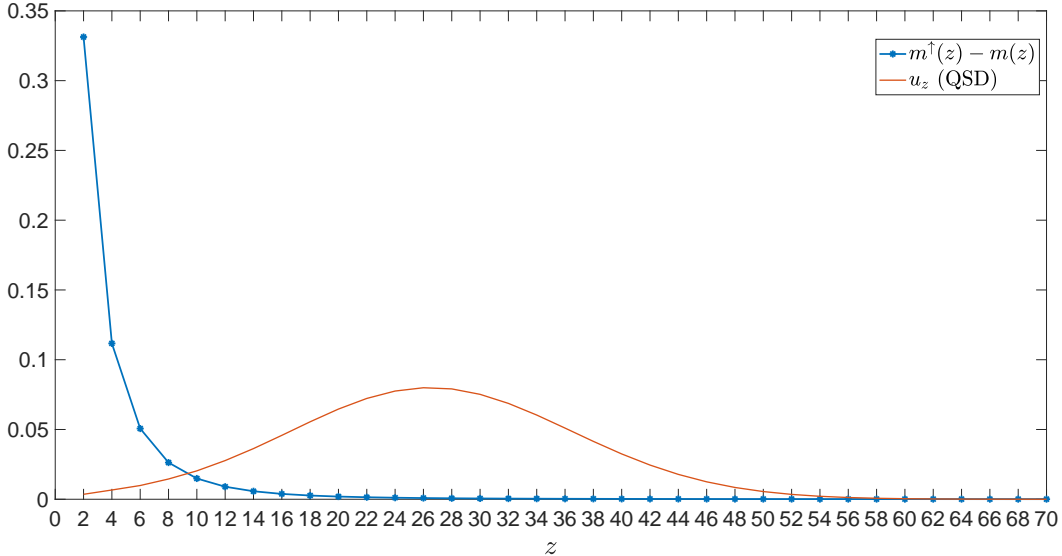}
\caption{\label{rick1} Absolute difference between the functions $m(z)$ and $m^\uparrow(z)$ plotted against the quasi-stationary distribution of the $Q$-process in the binary splitting Ricker model \eqref{ricker} with $r=1.2$ and carrying capacity $K=30$.}
\end{figure}

Let us return to the numerical example considered in Section \ref{sec:mle-pk} (binary splitting Ricker model with $r=1.2$ and carrying capacity $K=30$). In Figure \ref{rick3a} we plot the empirical distribution of $\hat{m}_{n}(z)$ for $z=28$ and $n=2000$ given $Z_n>0$, and the theoretical normal distribution implied by Theorem~\ref{thm:psdbp-Qcons and Qnorm}. In Figure \ref{rick3b} we plot the same quantities but in this case for  $z=8$. Comparing  Figures \ref{rick3a} and \ref{rick3b}, we see that when $z$ is closer to the absorbing state, the estimates and their limit are further away from the true value $m(z)$.
In Figure~\ref{rick1} we plot the difference $m^{\uparrow}(z)-m(z)$ as a function of $z$, together with the quasi-stationary distribution of the process. 
We see that $m^{\uparrow}(z)$ is substantially larger than $m(z)$ when $z$ is small. As indicated by the values of the quasi-stationary distribution, the population has a  small probability to be in these vulnerable states in the long run, therefore we can generally assume that $m^{\uparrow}(z)\approx m(z)$ for the values of $z$ observed in practice. However, we can legitimately question whether $Q$-consistency is a satisfactory property in other settings. We explore this matter in the next section.

\subsection{$C$-consistency versus $Q$-consistency --- the subcritical GW case} \label{sec:Ccons}


Recall that an estimator $\hat{\theta}_n$ is called $C$-consistent for a quantity $\theta$  in $\{Z_n\}$ when, for any $\varepsilon>0$, 
\begin{eqnarray*}\lim_{n\rightarrow\infty} \mbP(|\hat{\theta}_n-\theta|>\varepsilon\,|\,Z_n>0)=0;\end{eqnarray*}this is in contrast with $Q$-consistency, when
\begin{eqnarray*} \lim_{n\rightarrow\infty} \mbP(|\hat{\theta}_n-\theta^\uparrow|>\varepsilon\,|\,Z_n>0)=0,\end{eqnarray*}where $\theta^\uparrow$ is the quantity equivalent to $\theta$ in the $Q$-process $\{Z_n^\uparrow\}$ associated with $\{Z_n\}$.  In practice, we often study populations \emph{because} they are still alive, in which case  observations of the population should then be viewed as being generated under the condition $\{Z_n>0\}$. The resulting bias may be mitigated by a $C$-consistent estimator; however, in larger populations that do not face immediate danger, a $Q$-consistent estimator may be preferable.


To investigate when 
$C$-consistency is preferable to $Q$-consistency, we first require a $C$-consistent estimator for a branching process that experiences almost sure extinction. For this, 
we  focus on the GW case.  Recall that the MLE  for the mean offspring $m$ of a GW process based on the sample $\mathcal{Z}_n$, given by 
$$\hat{m}_n= \dfrac{\sum_{i=1}^{n} Z_i}{\sum_{i=1}^{n}Z_{i-1}}, $$
 is $C$-consistent in the supercritical case $m>1$ (see \cite[Theorem 7.2]{harris48}). 
 Pakes \cite{Pakes-1975} points out that almost nothing is known about the properties of $\hat{m}_n$ in the subcritical case $m< 1$, and that the elucidation of these properties would be of considerable interest. In his recent survey \cite{yanev08}, Yanev adds that it is well known that consistent and asymptotically normal estimators for parameters of GW processes exist only in the
supercritical case on the explosion set; otherwise, statistical inference in the subcritical case relies on the assumption of
an increasing initial population size or the presence of immigration.

 In the next theorem we show that, regardless the value of $m<1$, conditional on $Z_n>0$, $\hat{m}_n$ converges to 1 as $n\rightarrow\infty$. This confirms that $\hat{m}_n$ is not a $C$-consistent estimator, however it is $Q$-consistent (Corollary \ref{cor:GW-Qcons}).
\begin{thm}\label{thm:GW-Qcons and Qnorm}
If $m<1$, then
for any initial distribution $\vc\pi$ with finite second moment, the MLE $\hat{m}_n$ for $m$ in a GW process satisfies, for any $q<1$ and any $\varepsilon>0$,
\begin{equation}\label{GWm}\lim_{n\rightarrow\infty} \mbP_{\pi}(n^q\,|\hat{m}_n-1|>\varepsilon\,|\,Z_n>0)=0.\end{equation} 
\end{thm}

\begin{cor}\label{cor:GW-Qcons}
If $m<1$, then $\hat{m}_{n}$ is a $Q$-consistent estimator for $m$.
\end{cor}

\begin{figure}[t]
\centering\includegraphics[width=1\textwidth]{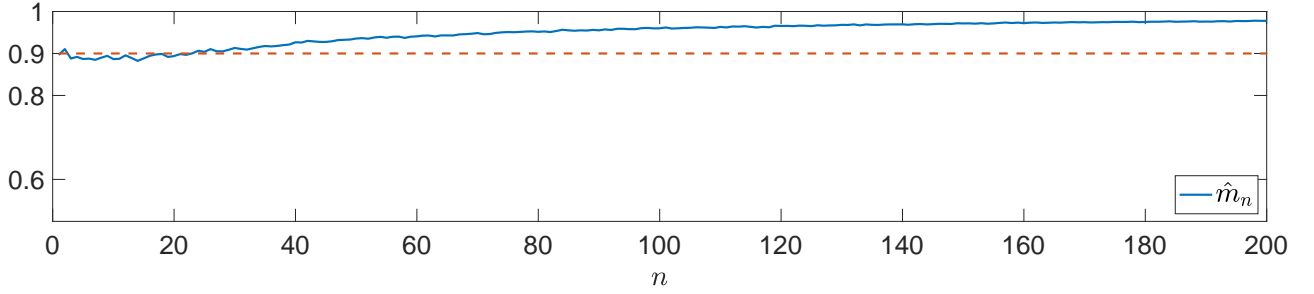}
\caption{\label{mhatGW}Successive values of $\hat{m}_{n}$ in a subcritical GW process with geometric offspring distribution with mean $m = 0.9$ and initial population size $Z_0 = 100$.}
\end{figure}
\medskip

Consider a GW process where the offspring distribution is geometric with mean $m = 0.9$ and the initial population size is $Z_0 = 100$. In Figure \ref{mhatGW}, we show successive values of $\hat{m}_{n}$ as $n$ increases, based on a simulated non-extinct trajectory.  Since $Z_0$ is large, for small values of $n$ the quality of the estimates benefits from the law of large numbers, but as $n$ increases, the population size becomes much smaller, and the condition $Z_n>0$ introduces a bigger bias in the estimates which start deviating from the true value. In the subcritical GW case, $Q$-consistency is therefore not always a desirable property. However, by exploiting the specific properties of the $Q$-process associated with a subcritical GW-process, we are able to construct  $C$-consistent estimators for $m$, as we describe now.

The $Q$-process $\{Z^{\uparrow}_n\}$ associated with a subcritical GW process $\{Z_n\}$ corresponds to the process of the generation sizes in the \emph{size-biased GW tree}; see for instance \cite{lyons}. If $\xi$ denotes the offspring distribution in $\{Z_n\}$, then the size-biased tree evolves as follows. 
At each generation, exactly one individual (the ``marked'' individual) reproduces according to the size-biased distribution of $\xi$, which is labelled SB$(\xi)$ and defined by $\mbP[\text{SB}(\xi)=k]=k\,p_k/m,$ $k\geq 0$, while the other individuals reproduce according to the offspring distribution $\xi$, all independently of each other (a more detailed construction is given in \cite{lyons}). Note that $Z^{\uparrow}_n\stackrel{d}{=}\textrm{SB}(Z_n)$.
Since $ \mbE(\text{SB}(\xi))=\mbE(\xi^2)/m$, we have
\begin{equation}\label{eq-zh}\mbE(Z^{\uparrow}_{\ell}\,|\,Z^{\uparrow}_{\ell-1})\;=\;(Z^{\uparrow}_{\ell-1}-1) \,m+\dfrac{\mbE(\xi^2)}{m}\;=\;Z^{\uparrow}_{\ell-1}\, m+\dfrac{\mbV(\xi)}{m},\quad \ell\geq1.\end{equation}


We first propose a $C$-consistent estimator for $m$ which holds for a class of offspring distributions $\xi$ satisfying the following assumption:
\begin{itemize}
\item[(A4)] There  exist constants $a$ and $b$ such that $$ \dfrac{\mbV(\xi)}{m} =a m+b\qquad  \textrm{with $a+b>0$}.$$
\end{itemize}
Observe that under Assumption (A4), by Equation \eqref{eq-zh} we have
$$m = \dfrac{\mbE(Z^{\uparrow}_{\ell}\,|\,Z^{\uparrow}_{\ell-1})-{b}}{Z^{\uparrow}_{\ell-1}+{a}}.$$ This suggests a way of modifying the observed data $Z_0,Z_1,\ldots,Z_n>0$ (which can be thought of as the generation sizes in the size-biased tree) so as to remove the bias introduced by the marked individuals. The next proposition formalises this concept.
%
\begin{prop}\label{Ccons1}If Assumptions (A3) and (A4) hold, then 
\begin{equation}\label{est1}
\tilde{m}_n:= \dfrac{\sum_{i=1}^{n} (Z_i-b)}{\sum_{i=1}^{n}(Z_{i-1}+a)}\end{equation}is a $C$-consistent estimator for $m$. In addition, for any $i\geq 1$ and $x\in \mathbb{R}$, we have
\begin{equation}\label{m_Cnorm1}\lim_{n\rightarrow\infty}\mbP_i\left(\{n/\tilde{\nu}^{2}\}^{1/2}\big(\tilde{m}_n-m\big)\leq x\,|\,Z_n>0\right)=\Phi(x),\end{equation}
where $\Phi(x)$ is the distribution function of a standard normal random variable, and $\tilde{\nu}^{2}$ is given in \eqref{nut2}.
\end{prop}
To the best of our knowledge, $\tilde{m}_n$ is the first $C$-consistent estimator for all values of $m$ (i.e., for the subcritical, critical and supercritical cases).
Examples of distributions satisfying Assumption (A4) are the Poisson, the geometric, and the Bernouilli distributions. In these cases, the estimator proposed in \eqref{est1} has a probabilistic interpretation: the constants $a$ and $b$ in \eqref{est1} remove the effect of the marked individuals in the size-biased tree. More precisely, if for each pair $(Z_i,Z_{i+1})$, $0\leq i\leq n-1$, we call $Z_i$ the progenitors and $Z_{i+1}$ the progeny, then
\begin{itemize}\item when $\xi\sim \text{Poi}(m)$, $a=0$ and $b=1$, and $\textrm{SB}(\xi)\stackrel{d}{=}1+\xi$: for each $i$, one progeny needs to be removed;
\item  when $\xi\sim \text{Geom}(1/(m+1))$, $a=1$ and $b=1$, and $\textrm{SB}(\xi)\stackrel{d}{=}1+\xi+\xi'$, where $\xi'$ is an independent copy of $\xi$:  for each $i$, one progenitor needs to be added and one progeny needs to be removed;
\item when $\xi\sim 2\, \text{Ber}(m/2)$, $a=-1$ and $b=2$, and $\textrm{SB}(\xi)=2$: for each $i$, one progenitor and two progeny need to be removed.
\end{itemize}

Another approach, which is not restricted to any particular class of offspring distributions, consists in interpreting the process $\{Z^{\uparrow}_n -1\}$ as a GW process with immigration, where the offspring distribution is $\xi$ and the immigration distribution is  SB$(\xi)-1$. We can then apply results on parameter estimation for branching processes with immigration, for example \cite{heyde1972estimation,heyde1972estimation_corr}. This leads to a second $C$-consistent estimator for $m$, as well as to a $C$-consistent estimator for the variance $\sigma^2= \mbV(\xi)$ of the offspring distribution.

\begin{prop}\label{Ccons2}If $m<1$, $p_1>0$, and Assumption (A3) holds, then 
\begin{equation}\label{est2}
\bar{m}_n:= 1-\dfrac{1}{2}\dfrac{\sum_{i=1}^{n} (Z_i-Z_{i-1})^2}{\sum_{i=1}^{n}(Z_{i-1}-\bar{Z}_n)^2},\qquad\text{where } \bar{Z}_n:=n^{-1}\,{\sum_{i=1}^n Z_{i-1}}\end{equation} is a $C$-consistent estimator for $m$.
In addition, for any  $i\geq 1$ and $x\in \mathbb{R}$, $\bar{m}_n$ satisfies
\begin{equation}\label{m_Cnorm2}\lim_{n\rightarrow\infty}\mbP_i\left(\{n/\bar{\nu}^{2}\}^{1/2}\big(\bar{m}_n-m\big)\leq x\,|\,Z_n>0\right)=\Phi(x),\end{equation}
where $\Phi(x)$ is the distribution function of a standard normal random variable, and $\bar{\nu}^{2}$ 
is given in \eqref{nu}.
A
$C$-consistent estimator for $\sigma^2$ is 
\begin{equation}\bar{\sigma}^2_n:=\bar{m}_n(1-\bar{m}_n)\bar{Z}_n.
\end{equation}
\end{prop}
Note that the asymptotic normality of $\bar{\sigma}^2_n$ is more challenging to obtain from \cite{heyde1972estimation} and is not pursued here. When $m\geq 1$, $\bar{m}_n$ is also a $C$-consistent estimator for $m$ but $\bar{\sigma}^2_n$ is not a $C$-consistent estimator for $\sigma^2$ (this follows from \cite[p. 1759]{wei90}).
Finally, we point out that it is likely that Assumption~(A3) can be weakened in Propositions \ref{Ccons1} and \ref{Ccons2}. 

\medskip

\begin{figure}[t]
\centering\includegraphics[width=0.8\textwidth]{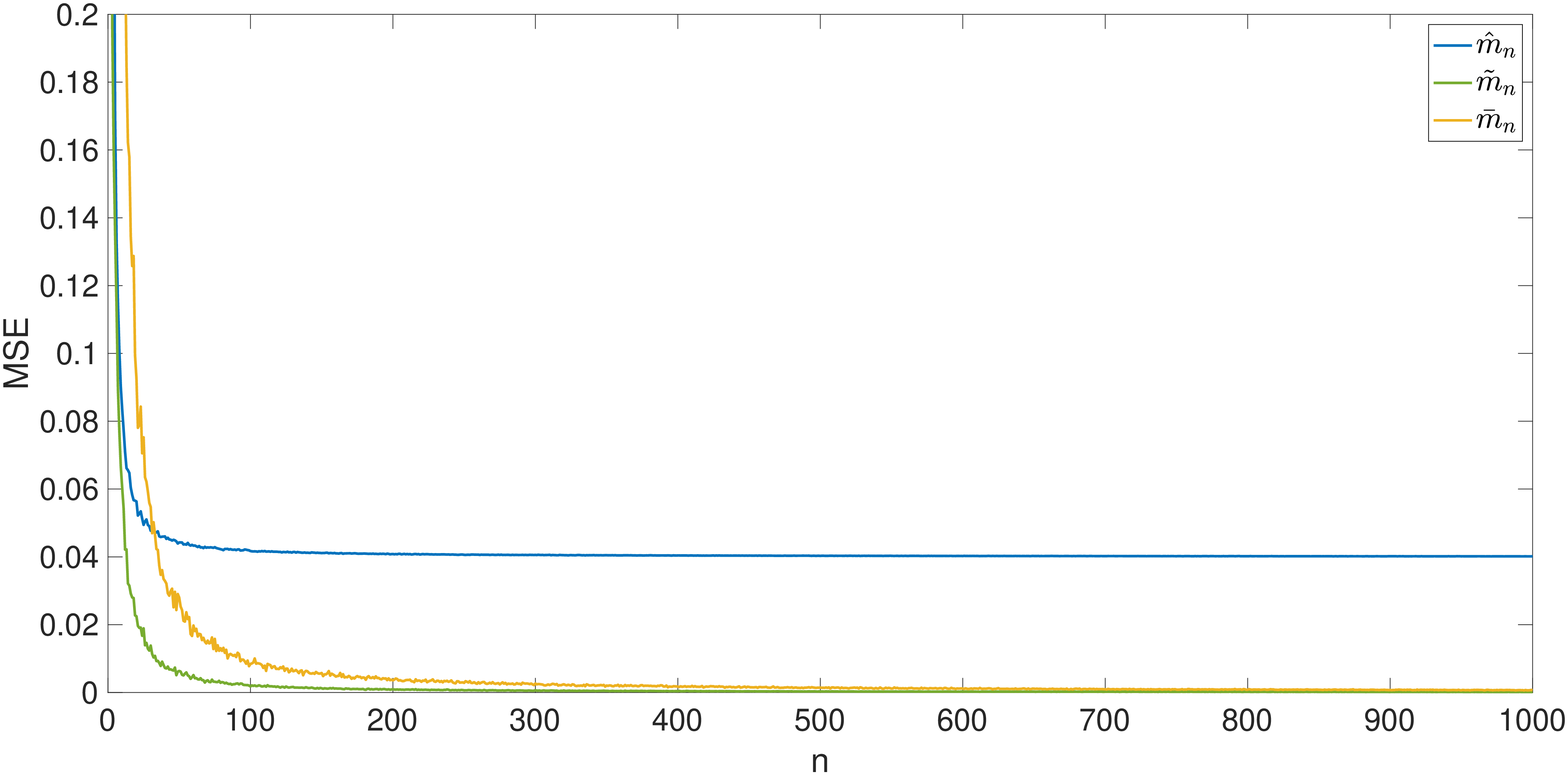}
\includegraphics[width=0.8\textwidth]{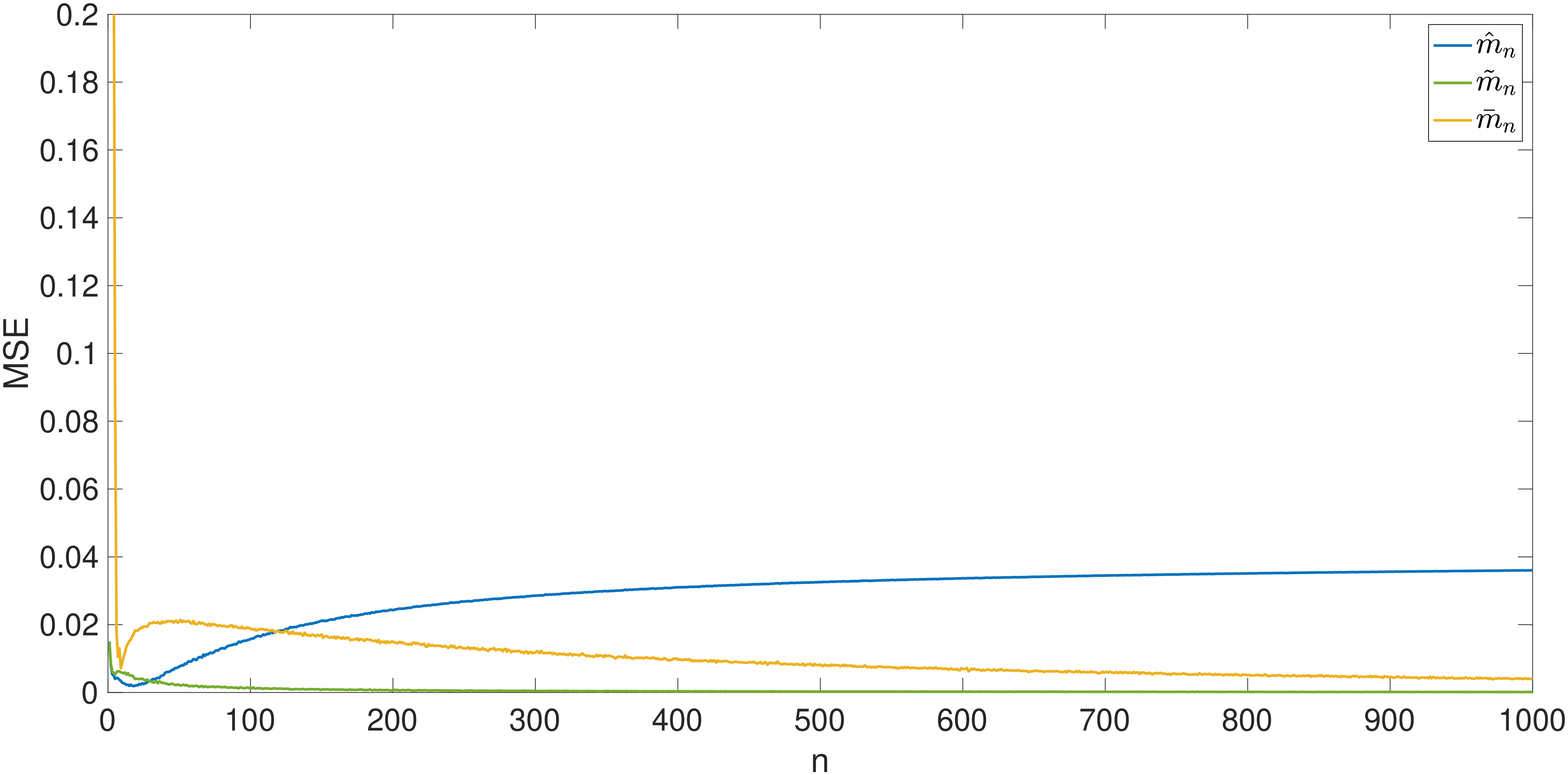}
\caption{\label{geom1MSE}Mean square error (MSE) of $\hat{m}_n$, $\tilde{m}_n$, and $\bar{m}_n$ computed from 500 trajectories of a GW process with geometric offspring distribution with mean $m=0.8$, with $Z_0=1$ (top) and $Z_0=100$ (bottom).}
\end{figure}

We illustrate the performance of the $C$-consistent estimators $\tilde{m}_n$ and $\bar{m}_n$ on a GW process where the offspring distribution is geometric with mean $m=0.8$. We consider trajectories of this process, conditional on $Z_n >0$ for $n=1,\ldots,1000$. These trajectories were simulated using the efficient method described in Section \ref{sec:coupling} below.
In Figure \ref{geom1} we plot the mean values of $\hat{m}_n$, $\tilde{m}_n$, and $\bar{m}_n$ averaged over 500 independent trajectories, with $Z_0=1$ and $Z_0=100$. Figure \ref{geom1MSE} shows the corresponding mean square errors (MSE).
Observe that if we start with a large population size, the $Q$-consistent estimator $\hat{m}_n$ has a smaller MSE for small values of $n$ ($n\leq28$), while the $C$-consistent estimators have a smaller MSE for larger values of $n$. 
In addition, we observe empirically that $\tilde{m}_n$ has a smaller MSE than $\bar{m}_n$. 

In conclusion, if a large population is observed for a small number of generations, the $Q$-consistent MLE $\hat{m}_n$ is preferable to the $C$-consistent estimators.  In contrast, if a non-extinct population is observed for a large number of generations, the sample becomes increasingly more biased by the condition $Z_n>0$, and the $C$-consistent estimators are more suitable.

\begin{figure}[t]
\centering\includegraphics[width=0.8\textwidth]{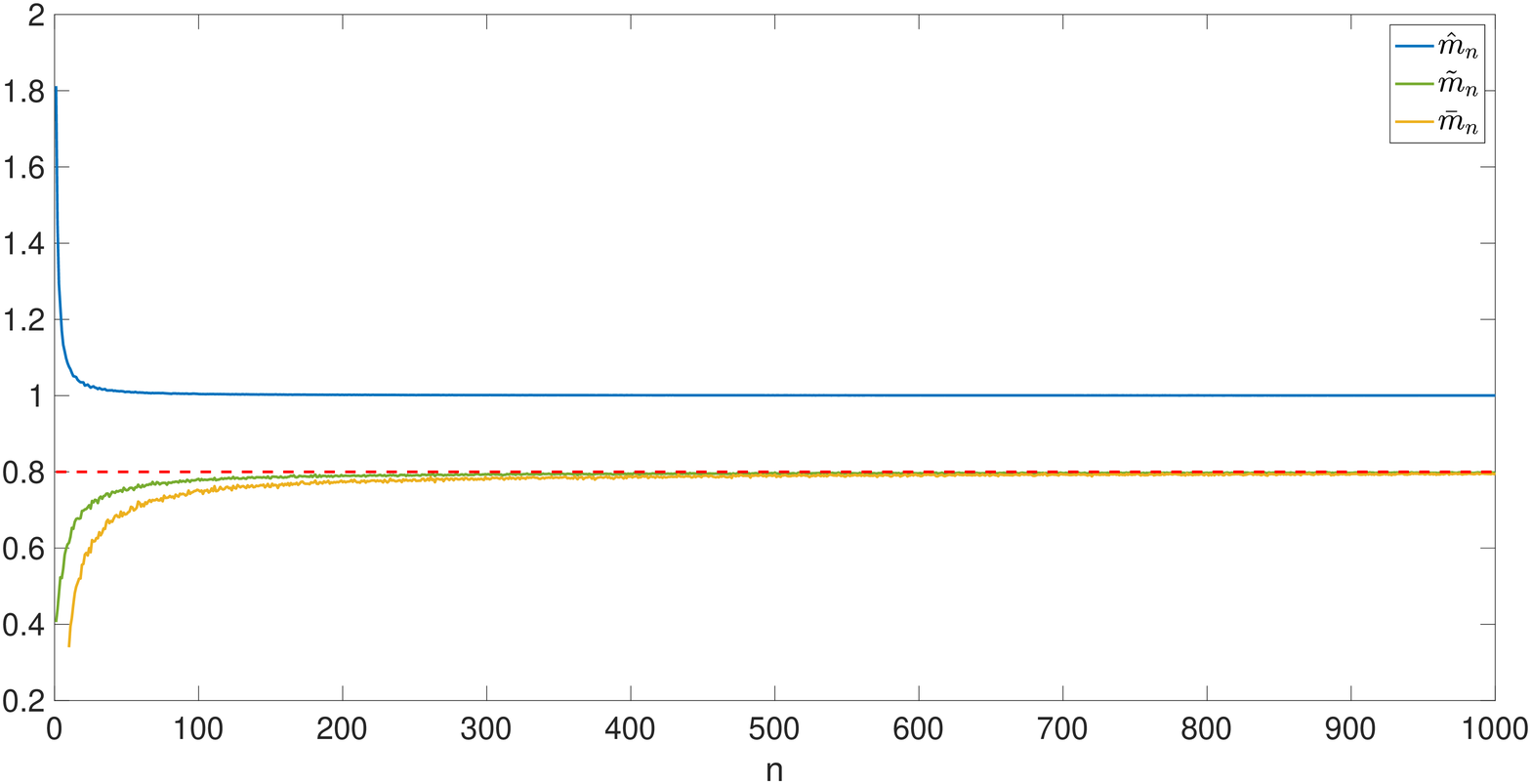}
\includegraphics[width=0.8\textwidth]{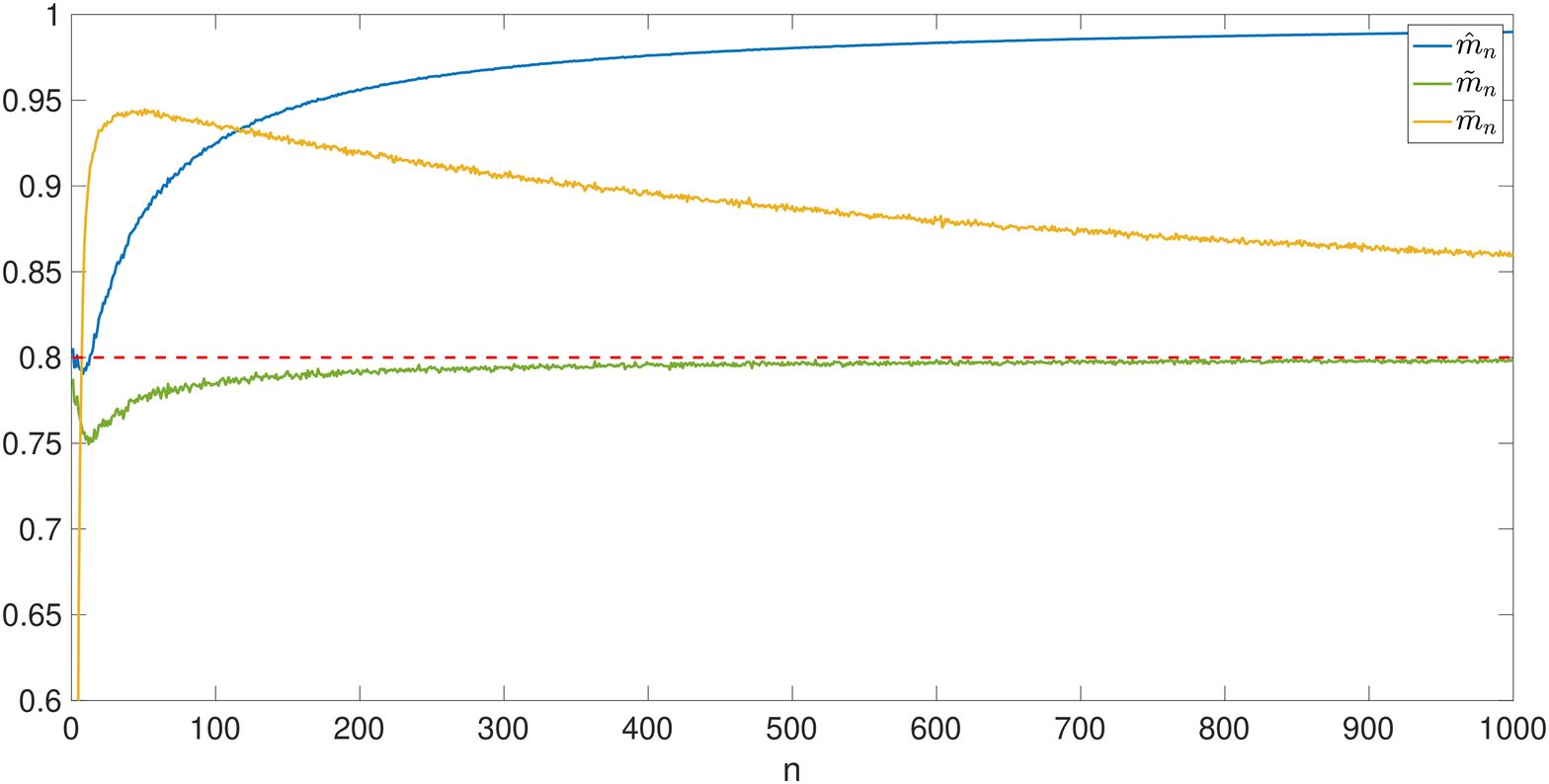}
\caption{\label{geom1} Mean values of $\hat{m}_n$, $\tilde{m}_n$, and $\bar{m}_n$ averaged over 500 trajectories of a GW process with geometric offspring distribution with mean $m=0.8$, with $Z_0=1$ (top) and $Z_0=100$ (bottom).}
\end{figure}

\subsection{Couplings of the $Q$-process}\label{sec:coupling}

One difficultly in analysing a branching process $\{Z_\ell\}_{0\leq \ell\leq n}$ conditional on $Z_n>0$ is that it evolves as  a time-\textit{inhomogeneous} Markov chain; we denote this time-{inhomogeneous} Markov chain by $\{Z_\ell^{(n)}\}_{0\leq \ell\leq n}$ (see Section \ref{subsec:Q-process}). 
To analyse the conditional asymptotic behaviour of the estimators in Sections \ref{sec:ass_prop} and \ref{sec:Ccons}, we need to analyse $\{Z_\ell^{(n)}\}$ for large values of $n$. It is however much simpler to manipulate the time-\textit{homogeneous} $Q$-process $\{Z^\uparrow_\ell\}$.
Here we show that there exists a coupling of $\{Z_\ell^{(n)}\}$ and $\{Z^\uparrow_\ell\}$ such that, for large $n$, the sample paths of the two processes coincide for a long time. We formalise this in Theorem \ref{coupling1}, after introducing  the necessary definitions.



For any $0\leq\ell\leq n $ and $\bm{x} \in \mbN^{\ell+1}$, we let
\begin{align*}
p^{(\ell,n)}(\bm{x})&:=\mbP(Z_0=x_0, Z_1=x_1,\ldots,Z_\ell=x_\ell\,  | Z_n >0)\\&=\mbP(Z_0^{(n)}=x_0, Z_1^{(n)}=x_1,\ldots,Z_\ell^{(n)}=x_\ell),\end{align*}
and \begin{align*}
p^{(\ell,\uparrow)}(\bm{x})&:=\lim_{n\rightarrow\infty}p^{(\ell,n)}(\bm{x})=\mbP(Z_0^\uparrow=x_0, Z_1^\uparrow=x_1,\ldots,Z_\ell^\uparrow=x_\ell).\\
\end{align*}
For each $n\geq 0$, a coupling of $\{Z_\ell^{(n)}\}$ and $\{{Z}^\uparrow_\ell \}$ is a random process $\{(\widehat{Z}_\ell^{(n)},\widehat{Z}^\uparrow_\ell )\}$  with associated probability measure $\widehat{\mbP}^{(n,\uparrow)}$, such that if for any $\ell\geq 0$ and $(\bm{s},\bm{t})\in (\mbN^{\ell+1})^2$, $$\widehat{p}^{(\ell,n,\uparrow)}(\bm{s},\bm{t}):=\widehat{\mbP}^{(n,\uparrow)}\left((\widehat{Z}_i^{(n)},\widehat{Z}^\uparrow_i )=(s_i,t_i), 0\leq i\leq \ell\right)$$
denotes the probability of the trajectory $(\bm{s},\bm{t})$ until time $\ell$,
then 
\begin{equation} \label{e1}
\sum_{\bm{t} \in \mbN^{\ell+1}} \widehat{p}^{(\ell,n,\uparrow)}(\bm{s},\bm{t}) = p^{(\ell,n)}(\bm{s}),\quad\textrm{and}\quad
\sum_{\bm{s}\in \mbN^{\ell+1}} \widehat{p}^{(\ell,n,\uparrow)}(\bm{s},\bm{t}) = p^{(\ell,\uparrow)}(\bm{t}),\qquad\textrm{for each $0\leq \ell \leq n$},
\end{equation}
(the marginal distributions are maintained) and, probabilities at consecutive times are related by
\begin{equation*}
\sum_{a\in \mbN} \sum_{b\in \mbN} p^{(\ell+1,n,\uparrow)}((\bm{s};a), (\bm{t};b)) = p^{(\ell,n,\uparrow)}(\bm{s},\bm{t}),
\end{equation*}where $(\bm{s};a):=(s_0,\ldots,s_\ell,a)$
(the probabilities are `inherited'). 
We define the uncoupling time as
\begin{equation}\label{uncoupling}\tau_n:=\min\{\ell\leq n: \widehat{Z}^{(n)}_\ell\neq \widehat{Z}^{\uparrow}_\ell\}.\end{equation}
If uncoupling has not occurred by generation $n$, we let $\tau_n=\infty$.
The next theorem states that there exists a sequence of couplings  such that $\tau_n$ is close to $n$ or equal to $\infty$ with high probability. These particular couplings are referred to as \emph{MEXIT couplings}, that is, maximal exit time couplings \cite{ernst2019mexit} (see also Section \ref{MEXIT} for more details). 
For $i\in\mbN$ we let $\widehat{\mbP}_i^{(n,\uparrow)}(\cdot)$ denote a MEXIT coupling of $\{ Z_\ell ^{(n)}| Z_0^{(n)}=i \}_{1 \leq \ell \leq n}$ and $\{ Z^\uparrow_\ell | Z^\uparrow_0 = i \}_{1\leq \ell \leq n},$
and we recall the definition of $\rho, \vc u, \vc v$, and $S$ in \eqref{PF1} and \eqref{uv}.

\begin{thm}\label{coupling1}
Suppose Assumptions (A1)--(A3) hold. Then, for any initial state $i>0$, there exists a sequence of couplings $\left\{\widehat{\mbP}_i^{(n,\uparrow)}(\cdot)\right\}_{n\geq 0}$ such that 
\begin{itemize}\item[(i)] for all $k \geq 1$,
\[
\lim_{n \to \infty} \widehat{\mbP}_i^{(n,\uparrow)}(\tau_n \leq n - k) = \frac{\rho^{-k}}{2} \sum_{j=1}^\infty u_j|\bm{e}_j^\top S^k \bm{1}|,
\]
\item[(ii)]
\[
\lim_{n \to \infty} \widehat{\mbP}_i^{(n,\uparrow)}(\tau_n < \infty) = \frac{1}{2}\sum_{j=1}^\infty u_j |1-v_j|,
\]
and, 
\item[(iii)] for any $q>0$, there exist constants $C(i,q)$ and $N(i,q)$ such that 
\[
\widehat{\mbP}_i^{(n,\uparrow)}\left(\tau_n \leq n - C(i,q)\log n \right) \leq 1/n^q
\]
for all $n \geq N(i,q)$.
\end{itemize}
\end{thm}
We use Theorem~\ref{coupling1}~(iii) to prove the results  in Sections  \ref{sec:ass_prop} and \ref{sec:Ccons}. In addition, Theorem~\ref{coupling1}~(i) also suggests an efficient method of approximate simulation of trajectories of $\{Z_\ell\}$ conditional on $Z_n>0$ (i.e. trajectories of $\{Z_\ell^{(n)}\}$) for large $n$. 
The idea is that it is significantly less computationally intensive to simulate trajectories of the (time-homogeneous) $Q$-process $\{Z^\uparrow_\ell\}$ than it is to simulate trajectories of the (time-inhomogeneous) process $\{Z_\ell^{(n)}\}$. 
According to Theorem~\ref{coupling1}~(i), to obtain a non-extinct trajectory of (large) length $n$, we can choose $k$ such that $d(k):=({\rho^{-k}}/{2}) \sum_{j\geq1} u_j|\bm{e}_j^\top S^k \bm{1}|$ is arbitrarily close to zero, simulate the $Q$-process up to generation $n-k$, and then proceed by simulating the last $k$ generations of the original process conditional on $Z_n>0$  using another method ---such as the multilevel splitting method (see Appendix A). The controlled error $d(k)$ is the total variation distance between the distribution of the simulated trajectory and the distribution of a trajectory of the original process conditional on $Z_n>0$. 

In Section \ref{sec:Ccons}, we used this method to simulate trajectories of length up to $n=1000$ from a GW process with geometric offspring distribution ($m=0.8$). In that particular example, by choosing $k=50$, the error is $d(k)=3.4481\cdot 10^{-6}$, and by choosing $k=100$, the error reduces to $d(k)=2.4607\cdot 10^{-11}$. This method of approximate simulation is very accurate and also significantly increases computational efficiency.

\vspace{0.5cm}



\section{Proofs}\label{proofs}
\subsection{Proofs of the results in Section \ref{sec:mle-pk}}\label{ape:mle-pk}
In this section, we prove Propositions \ref{prop:mle-m} and \ref{prop:moments-mle}. The result in Proposition \ref{prop:mle-m} is part of  Corollary~\ref{cor:mle-m-sigma} below. More precisely, we derive the MLE for the offspring distribution at population size $z$ in a PSDBP,
$\boldsymbol{p}(z)$ (Theorem \ref{thm:mle-pk}), as well as for its mean $m(z)$ and variance $\sigma^2(z)$ (Corollary~\ref{cor:mle-m-sigma}). These estimators are based on the observation of the entire family tree up to some generation. Specifically, we consider the sample $\mathcal{Z}_n^*=\{Z_i(k):i=0,\ldots,n-1,\ k\in\N_0\}$, where $Z_i(k)$ represents the number of individuals at generation $i$ that have exactly $k$ offspring, that is,
$$Z_i(k)=\sum_{j=1}^{Z_i}\ind{\xi_{ij}(Z_i)=k},\quad i=0,\ldots,n-1,\ k\in\N_0.$$ Observe that
$$\sum_{k=0}^\infty Z_i(k)=Z_i,\quad \sum_{k=0}^\infty kZ_i(k)=Z_{i+1},\qquad i=0,\ldots,n-1.$$
Thus, because $Z_i$ is finite for each $i\in\{0,\ldots,n\}$, only a finite number of elements in the sequence $(Z_i(k))_{k\in\N_0}$ are non-null. 
\begin{thm}\label{thm:mle-pk}
The likelihood function based on the sample $\mathcal{Z}_n^*$ is
\begin{equation*}
\mathcal{L}(\mathcal{Z}_n^*)=\prod_{i=0}^{n-1}\frac{Z_i!}{\prod_{k=0}^\infty Z_i(k)!}\prod_{k=0}^\infty p_k(Z_i)^{Z_i(k)},
\end{equation*}
and for every $z\in\N$ and $k\in\N_0$, the MLE of $p_k(z)$ on the set $\{j_n(z)>0\}$ is  
\begin{equation}\label{eq:mle-pk}
\hat{p}_{k,n}(z)=\frac{\sum_{i\in J_{n}(z)} Z_{i}(k)}{j_n(z)z}=\frac{\sum_{i=0}^{n-1} Z_{i}(k)\ind{Z_i=z}}{j_n(z)z},
\end{equation}
where $J_n(z)=\{i=0,\ldots,n-1:Z_i=z\}$ is the set of the generation indexes when the process visits the state $z$.
\end{thm}

\begin{pro}
Using the Markov property, for the likelihood function we have
\begin{align*}
\mbP(Z_i(k)=z_i(k),&\ i=0,\ldots,n-1,\ k\in\N_0)=\\
&=\prod_{i=0}^{n-1} \mbP(Z_i(k)=z_i(k),\  k\in\N_0 | \ Z_{h}(k) = z_{h}(k),\ h=0,\ldots i-1, k\in\N_0 )\\
&= \prod_{i=0}^{n-1} \mbP(Z_i(k)=z_i(k),\  k\in\N_0 | \ Z_i=z_i,\ Z_{i-1}(k) = z_{i-1}(k), k\in\N_0)\\
&= \prod_{i=0}^{n-1} \frac{z_i!}{\prod_{k=0}^\infty z_i(k)!}\cdot \prod_{k=0}^\infty p_k(z_i)^{z_i(k)},
\end{align*}
where $z_i(k)\in\N_0$, for each $i=0,\ldots,n-1$, $k\in\N_0$ and $z_i=\sum_{k=0}^\infty k z_{i-1}(k)$.

For the second part, conditionally on $\{Z_i(k)=z_i(k):i=0,\ldots,n-1,\ k\in\N_0\}$, there exist $i_1,\ldots,i_{L(n)}\in\N$ such that $z_{i_1},\ldots,z_{i_{L(n)}}$ are all the population sizes observed in the sample $z_0,\ldots,z_{n-1}$ but satisfying $z_{i_j}\neq z_{i_l}$ if $j\neq l$, and $1\leq L(n)\leq n$. Thus, with this notation the MLEs of $p_{k}(z_i)$, denoted $\hat{p}_{k,n}(z_i)$, $k\in\N_0$, and $\hat{\boldsymbol{p}}_{n}(z_i)=(\hat{p}_{k,n}(z_i))_{k\in\N_0}$, $i=0,\ldots,n-1$, are given by
\begin{equation*}
\left\{\hat{\boldsymbol{p}}_{n}(z_i)\right\}_{i=0,\ldots,n-1}=\arg\max_{\bar{\boldsymbol{p}}(z_i):i=0,\ldots,n-1} f_n(\bar{\boldsymbol{p}}(z_0),\ldots,\bar{\boldsymbol{p}}(z_{n-1})),
\end{equation*}
with
\begin{equation*}
f_n(\bar{\boldsymbol{p}}(z_0),\ldots,\bar{\boldsymbol{p}}(z_{n-1}))=\sum_{i=0}^{n-1}\sum_{k=0}^\infty z_{i}(k)\log(\bar{p}_k(z_{i}))=\sum_{l=1}^{L(n)}\sum_{j\in  J_{n}(z_{i_l})} \sum_{k=0}^\infty z_{j}(k)\log(\bar{p}_k(z_{j})),
\end{equation*}
and $\bar{\boldsymbol{p}}(z_i)$ are taken to be probability distributions on $\N_0$, for $i=0,\ldots,n-1$. Now, it is easy to check that the values in \eqref{eq:mle-pk} are the solutions of the system of equations
\begin{align*}
\frac{\partial}{\partial p_h(z_{i_l})}\left(f_n(\bar{\boldsymbol{p}}(z_0),\ldots,\bar{\boldsymbol{p}}(z_{n-1})) +\sum_{l=1}^{L(n)} \lambda_{i_l} \left(1- \sum_{k=0}^\infty \bar{p}_k(z_{i_l})\right)\right) &= 0, \\
\frac{\partial}{\partial\lambda_{i_l}}\left(f_n(\bar{\boldsymbol{p}}(z_0),\ldots,\bar{\boldsymbol{p}}(z_{n-1})) +\sum_{l=1}^{L(n)} \lambda_{i_l} \left(1- \sum_{k=0}^\infty \bar{p}_k(z_{i_l})\right)\right) &= 0,
\end{align*}
for $h\in\N_0$ and $l=1,\ldots,L(n)$. Finally, note that to prove that the values in \eqref{eq:mle-pk} are a maximum of the function $f_n(\bar{\boldsymbol{p}}(z_0),\ldots,\bar{\boldsymbol{p}}(z_{n-1}))$, it is equivalent to prove that they are a maximum of the function $y_{n-1}^{-1}f_n(\bar{\boldsymbol{p}}(z_0),\ldots,\bar{\boldsymbol{p}}(z_{n-1}))$, where $y_{n-1}=\sum_{i=0}^{n-1}z_i$. Since $x\in [0,\infty)\mapsto -\log(x)$ is a convex function, by applying Jensen's inequality we have that for any family of probability distributions $\bar{\boldsymbol{p}}(z_0),\ldots,\bar{\boldsymbol{p}}(z_{n-1})$,
\begin{align*}
y_{n-1}^{-1}&f_n(\bar{\boldsymbol{p}}(z_0),\ldots,\bar{\boldsymbol{p}}(z_{n-1}))-y_{n-1}^{-1}f_n(\boldsymbol{\hat{p}}_{n}(z_0),\ldots,\boldsymbol{\hat{p}}_{n}(z_{n-1}))=\\
&=\sum_{l=1}^{L(n)} \frac{j_{n}(z_{i_l})z_{i_l}}{y_{n-1}}\sum_{k=0}^\infty \log\bigg(\frac{\bar{p}_k(z_{i_l})}{\hat{p}_{k,n}(z_{i_l})}\bigg)\frac{\sum_{j\in J_{n}(z_{i_l})} z_{j}(k)}{j_{n}(z_{i_l})z_{i_l}}\\
&\leq\sum_{l=1}^{L(n)} q_n(z_{i_l}) \log\bigg(\sum_{k=0}^\infty\frac{\bar{p}_k(z_{i_l})}{\hat{p}_{k,n}(z_{i_l})}\cdot\hat{p}_{k,n}(z_{i_l})\bigg)=0,
\end{align*}
with $0<q_n(z_{i_l})=\frac{j_{n}(z_{i_l})z_{i_l}}{y_{n-1}}< 1$, for each $l=1,\ldots,L(n)$.
\end{pro}

\vspace{0.5cm}

\begin{cor}\label{cor:mle-m-sigma}
For every $z\in\N$, the MLEs of $m(z)$, and $\sigma^2(z)$ on the set $\{j_n(z)>0\}$ are 
\begin{equation*}\label{eq:mle-m-sigma}
\begin{aligned}
\hat{m}_{n}(z)&=\frac{\sum_{i\in J_{n}(z)} Z_{i+1}}{j_n(z)z}=\frac{\sum_{i=0}^{n-1} Z_{i+1}\ind{Z_i=z}}{j_n(z)z},\\
\hat{\sigma}_{n}^2(z)&=\sum_{k=0}^\infty (k-\hat{m}_{n}(z))^2 \hat{p}_{k,n}(z).
\end{aligned}
\end{equation*}
For every fixed $z\in\N$, $\hat{m}_n(z)$ is also the MLE on the set $\{j_n(z)>0\}$ based on the sample $\mathcal{Z}_n:=\{Z_0,\ldots,Z_{n}\}$. 
\end{cor}

\begin{pro}
The result follows from Theorem \ref{thm:mle-pk} and from the invariance of the MLEs under continuous transformations. For the offspring mean function, conditionally on $\{Z_i(k)=z_i(k):i=0,\ldots,n-1,\ k\in\N_0\}$,  note that
\begin{equation*}
\hat{m}_{n}(z_i)=\sum_{k=0}^\infty k \hat{p}_{k,n}(z_i)=\frac{\sum_{k=0}^\infty\sum_{l\in J_{n}(z_{i})}k z_l(k) }{j_n(z_i)z_{i}}=\frac{\sum_{l\in J_{n}(z_{i})} z_{l+1}}{j_n(z_i)z_{i}}.
\end{equation*}The fact that $\hat{m}_n(z)$ is also the MLE on the set $\{j_n(z)>0\}$ based on the sample $\mathcal{Z}_n$ can be proved by using Lemma 2.13.2 in \cite{jager} and the fact that $\hat{m}_n(z)$ is a measurable function of the variables $Z_0,\ldots,Z_n$.\label{rem:MLE-zn}
\end{pro}

\medskip
%
%
\begin{remark}\label{rem:binary}
In the case of a PSDBP with binary splitting reproduction ($p_0(z)+p_2(z)=1$), there is a one-to-one correspondence between the samples $\mathcal{Z}_n$ and $\mathcal{Z}_n^*$. Indeed, 
$$Z_k(0)=Z_k-\frac{Z_{k+1}}{2},\quad Z_k(2)=\frac{Z_{k+1}}{2},\quad Z_k(h)=0,\quad\text{ for each }k=0,\ldots,n-1;\ h\in\N\backslash\{2\}.$$
Thus, we can compute all MLEs in Theorem \ref{thm:mle-pk} and Corollary \ref{cor:mle-m-sigma} based on the sample $\mathcal{Z}_n$. Moreover, by the last part of Corollary \ref{cor:mle-m-sigma} and the invariance of the MLE, we conclude that the MLEs of $p_0(z)$, $p_2(z)$ and $\sigma^2(z)$ based on the sample $\mathcal{Z}_n$ on the set $\{j_n(z)>0\}$ are, respectively,
$$\hat{p}_{0,n}(z)=1-\frac{\hat{m}_n(z)}{2},\quad \hat{p}_{2,n}(z)=\frac{\hat{m}_n(z)}{2},\quad \text{ and }\quad\hat{\sigma}_n^2(z)=2\hat{m}_n(z)\left(1-\hat{m}_n(z)\right).$$ 
\end{remark}

\begin{Prf}[Proposition \ref{prop:moments-mle}.]
Given $z_0,\ldots,z_{n-1}\in\N$, we have
\begin{align*}
\mbE(\hat{m}_{n}(z)\ind{j_n(z)>0}&\,|\,Z_0=z_0,\ldots,Z_{n-1}=z_{n-1})=\\
&=\mbE\left(\frac{Z_{n}\ind{Z_{n-1}=z}}{z\sum_{i=0}^{n-1}\ind{Z_i=z}}\,\Big|\,Z_0=z_0,\ldots,Z_{n-1}=z_{n-1}\right)\\
&\phantom{=}+\mbE\left(\frac{\sum_{i=0}^{n-2}Z_{i+1}\ind{Z_i=z}}{z\sum_{i=0}^{n-1}\ind{Z_i=z}}\,\Big|\,Z_0=z_0,\ldots,Z_{n-1}=z_{n-1}\right)\\
&=\frac{z_{n-1}m(z_{n-1})\ind{z_{n-1}=z}}{z\sum_{i=0}^{n-1}\ind{z_i=z}}+\frac{\sum_{i=0}^{n-2}z_{i+1}\ind{z_i=z}}{z\sum_{i=0}^{n-1}\ind{z_i=z}},
\end{align*}
and consequently,
\begin{align*}
\mbE\left(\hat{m}_{n}(z)\,|\,j_n(z)>0\right)
&=\frac{1}{\mbP(j_n(z)>0)} \mbE\left(\mbE\left(\hat{m}_{n}(z)\ind{j_n(z)>0}\,|\,Z_0,\ldots,Z_{n-1}\right)\right)\\
&=\frac{1}{\mbP(j_n(z)>0)} \mbE\left(\frac{Z_{n-1}m(Z_{n-1})\ind{Z_{n-1}=z}}{z\sum_{i=0}^{n-1}\ind{Z_i=z}}+\frac{\sum_{i=0}^{n-2}Z_{i+1}\ind{Z_i=z}}{z\sum_{i=0}^{n-1}\ind{Z_i=z}}\right)\\
&= \mbE\left(\frac{Z_{n-1}m(Z_{n-1})\ind{Z_{n-1}=z}}{zj_n(z)}+\frac{\sum_{i=0}^{n-2}Z_{i+1}\ind{Z_i=z}}{zj_n(z)}\,\Big|\, j_n(z)>0\right)\\
&= m(z)\mbE\left(\frac{\ind{Z_{n-1}=z}}{j_n(z)}\,\Big|\, j_n(z)>0\right)+\mbE\left(\frac{\sum_{i=0}^{n-2}Z_{i+1}\ind{Z_i=z}}{zj_n(z)}\,\Big|\, j_n(z)>0\right).
\end{align*}
\end{Prf}

\subsection{Proofs of the results in Section \ref{sec:ass_prop}}\label{ape:ass_prop}

\subsubsection{Martingale central limit theorem}\label{mclt}
In the proofs of the results in Section \ref{sec:ass_prop}, we make use of the following version of the martingale central limit theorem (see for  instance \cite{hall}). Let $M_n$  be a  mean-$0$ $\mathbb{R}^d$-valued vector martingale in $L^2$ (that is, each coordinate forms a real-valued
martingale) with increments $Y_k=M_k-M_{k-1}$. Define the scaled random variables $\bar{M}_n=n^{-1/2} M_n$. Assume
\begin{equation}\label{clt-cond1}
\dfrac{1}{n}\sum_{k=1}^{n} \mathbb{E}(Y_k Y_k^\top\,|\,\mathcal{F}_{k-1})\longrightarrow \Gamma\qquad\text{in probability}\end{equation}for a symmetric, nonnegative definite $d\times d$ matrix $\Gamma$, and for any $\varepsilon>0$,
\begin{equation}\label{clt-cond2}
\dfrac{1}{n}\sum_{k=1}^{n} \mathbb{E}(|Y_k|^2 \,\ind{|Y_k|\geq \varepsilon \sqrt{n}}\,|\,\mathcal{F}_{k-1})\longrightarrow 0\qquad\text{in probability}.\end{equation}
Then $\bar{M}_n$ converges weakly to a $d$-dimensional centred normal random variable with covariance matrix $\Gamma$. 

\subsubsection{The proportion of time spent in a transient state}

Let $z\geq 1$ be a transient state of an absorbing Markov chain $\{Z_n\}_{n\in\N_0}$. 
 The next lemma uses Theorem  \ref{coupling1} to recover Gosselin's result \cite[Theorem 3.1(c)]{Gosselin-2001} on the proportion of time $\{Z_n\}$ spends in state $z$ up to time $n$, which is used to prove Theorem \ref{thm:psdbp-Qcons and Qnorm}.
 
 \begin{lem}\label{proof1}
 Under Assumptions (A1)--(A3),
 for all $i ,z\geq 1$ and $\varepsilon>0$, we have 
 \[
 \mbP_i \left( \left| \frac{\sum_{\ell=1}^n \ind{ Z_\ell = z}}{n} - u_z v_z \right| > \varepsilon 
 , \Big| \, Z_n >0\right) \to 0
 \]
 as $n \to \infty$.
 \end{lem}

 \begin{pro}
Recall from Section \ref{sec:coupling} that $\widehat{\mbP}_i^{(n,\uparrow)}(\cdot)$ denotes the probability measure associated with a MEXIT coupling  of $\{ Z_\ell ^{(n)}| Z_0^{(n)}=i \}_{1 \leq \ell \leq n}$ and $\{ Z^\uparrow_\ell | Z^\uparrow_0 = i \}_{1\leq \ell \leq n}$, and $\tau_n$ is the uncoupling time defined in \eqref{uncoupling}.
 Observe that, on the paths where $\tau_n > n - C(i,q) \log n$, 
\begin{equation}\label{diff_zero}
 n^{-1} \left| \sum_{\ell=1}^n \ind{ \widehat{Z}^{(n)}_\ell =z } - \sum_{\ell=1}^n \ind{ \widehat{Z}_\ell^{\uparrow} =z } \right| \leq \frac{C(i,q)\log n }{n} \to 0
\end{equation}
as $n \to \infty$, where $C(i,q)$ is defined in Theorem \ref{coupling1} (iii).
 We thus have, for any $\varepsilon>0$,
 \begin{align*}
\mbP_i &\left( \left| \frac{\sum_{\ell=1}^n \ind{ Z_\ell = z}}{n} - u_z v_z  \right| > \varepsilon 
 \, \Big| \, Z_n >0\right)=\mbP_i \left(\left| \frac{\sum_{\ell=1}^n \ind{ Z_\ell^{(n)} = z}}{n} - u_z v_z  \right| > \varepsilon 
 \right)\\
  &\leq \widehat{\mbP}_i^{(n,\uparrow)}\left( \left| \frac{\sum_{\ell=1}^n \ind{  \widehat Z^\uparrow_\ell = z}}{n} - u_z v_z  \right| > \frac{\varepsilon}{2} \right) + \widehat{\mbP}_i^{(n,\uparrow)} \left(  \left| \frac{\sum_{\ell=1}^n \ind{  \widehat Z_\ell^{(n)} = z}}{n} -\frac{\sum_{\ell=1}^n \ind{  \widehat Z^\uparrow_\ell = z}}{n}  \right| > \frac{\varepsilon}{2} \right)\\
 &\leq \widehat{\mbP}_i^{(n,\uparrow)}\left(  \left| \frac{\sum_{\ell=1}^n \ind{  \widehat Z^\uparrow_\ell = z}}{n} - u_z v_z \right| > \frac{\varepsilon}{2} \right) \\&+ \widehat{\mbP}_i^{(n,\uparrow)}\left(  \left| \frac{\sum_{\ell=1}^n \ind{  \widehat Z_\ell^{(n)} = z}}{n} -\frac{\sum_{\ell=1}^n \ind{  \widehat Z^\uparrow_\ell = z}}{n}  \right|   > \frac{\varepsilon}{2},\,\tau_n > n - C(i,q) \log n \right)\\&+\widehat{\mbP}_i^{(n,\uparrow)}\left(\tau_n \leq n - C(i,q)\log n\right). 
 \end{align*}As $n \to \infty$, the second term vanishes by \eqref{diff_zero}, and the third term vanishes by Theorem \ref{coupling1} (iii).
 The result then follows from the ergodic theorem for the positive recurrent Markov chain $\{ Z^\uparrow_\ell\}$.

%
%
%
%
%
%
%
%
 \end{pro}

\subsubsection{Proof of Theorem \ref{thm:psdbp-Qcons and Qnorm} and Corollary \ref{cor:psdbp-Qcons and Qnorm}}\label{pr_PSDBP}

In this subsection, we implicitly place ourselves in the probability space of a MEXIT coupling $\{(\widehat{Z}_\ell^{(n)},\widehat{Z}^\uparrow_\ell )\}$ defined in Section \ref{sec:coupling}, and we simply write $\mbP$ for the associated measure $\widehat{\mbP}^{(n,\uparrow)}(\cdot)$. We also assume that Assumptions (A1)--(A3) hold.

Combining Lemma \ref{proof1}, the continuous mapping theorem, and Slutsky's theorem, we see that, if $Z \sim N(0,z{\sigma^2}^{\uparrow}(z))$ where ${\sigma^2}^{\uparrow}(z)$ is defined in \eqref{sig2b}, then proving \eqref{mz_Qnorm} in Theorem \ref{thm:psdbp-Qcons and Qnorm} is equivalent to proving that, for any initial state $i\geq1$ and every $z\in\N$, conditional on $Z_0=i$, 
\begin{equation}\label{Ynn}
\widehat Y_n^{(n)}(z):= \frac{\sum_{\ell=1}^n Z^{(n)}_{\ell} \ind{ Z_{\ell-1}^{(n)}=z} -z m^{\uparrow}(z) \sum_{\ell=1}^n \ind{ Z_{\ell-1}^{(n)}=z}  }{(z\,u_z v_z n)^{1/2}} \stackrel{d}{\to} Z.
\end{equation}
We define the quantity similar to $\widehat Y_n^{(n)}(z)$ in the $Q$-process,
\[
\widehat Y_n^\uparrow(z):= \frac{\sum_{\ell=1}^n Z^\uparrow_{\ell} \ind{ Z_{\ell-1}^\uparrow=z} -z m^{\uparrow}(z) \sum_{\ell=1}^n \ind{ Z_{\ell-1}^\uparrow=z}  }{(z\,u_z v_z n)^{1/2}}.
\]We proceed by showing that $\widehat Y_n^{(n)}(z) - \widehat Y_n^\uparrow(z) \stackrel{P}{\to} 0$ and $\widehat Y_n^\uparrow(z) \stackrel{d}{\to} Z$ in the next two lemmas. 
 
\begin{lem}\label{le1}
For any initial state $i\geq1$ and every $z\in\N$, conditional on $Z_0=i$, we have $\hat Y_n^{(n)}(z) - \widehat Y_n^\uparrow(z) \stackrel{P}{\to} 0$.
\end{lem} 
\begin{pro}
Recall that $\tau_n$ is the uncoupling time defined in \eqref{uncoupling} and $C(i,q)$ is the constant defined in Theorem \ref{coupling1} (iii).
If the events 
\begin{align*}
A_{1,n} &:=\{ \tau_n > n - C(i,q) \log n \}, \\ 
A_{2,n} &:=\{ \nexists \ell: (n - C(i,q) \log n) \leq \ell \leq n,\,  Z^{\uparrow}_{\ell-1} = z, Z^{\uparrow}_{\ell} > n^{1/4} \},  \\
A_{3,n}&:=\{ \nexists \ell: (n - C(i,q) \log n) \leq \ell \leq n,\, Z^{(n)}_{\ell-1} = z, Z^{(n)}_{\ell} > n^{1/4} \} 
\end{align*}
all occur, then 
\[
|\hat Y_n^{(n)}(z) - \hat Y_n^\uparrow(z)| \leq \frac{C(i,q)  \log n \,(n^{1/4}+z m^{\uparrow}(z))}{(z\,u_z v_z n)^{1/2}},
\]
which becomes arbitrarily small for large $n$. 
Thus, if we can show that $\mbP(A_{1,n}), \mbP(A_{2,n}), \mbP(A_{3,n}) \to 1$ as $n \to \infty$ then the result is proved.

By Theorem \ref{coupling1} (iii) we have $\mbP(A_{1,n}) \geq 1 - n^{-q} \to 1$. 
To show that  $\mbP(A_{2,n}) \to 1$, let $X^\uparrow$ (resp. $X$) be a random variable with $\mbP(X^\uparrow = j) = Q^\uparrow_{zj}$ (resp. $\mbP(X = j) = Q_{zj}$). 
Using Markov's inequality in the first step, the fact that for any $\nu>\nu^*$ there exists $C$ such that  $v_j/j^\nu < C$ for all $j$ (Lemma \ref{Goss2}) in the second step, and Minkowski's inequality in the third step, we obtain  
\begin{align*}
\mbP(X^\uparrow > n^{1/4} ) \leq \frac{\mbE\left( (X^\uparrow)^r\right)}{n^{r/4}} \leq \frac{C \,\mbE\left(X^{r+\nu}\right)}{\rho\,v_z n^{r/4}} \leq \frac{C \,z^{r+\nu}\, \mbE\left(\xi(z)^{r+\nu}\right)}{\rho\,v_z n^{r/4}},
\end{align*}where $\xi(z)$ is the offspring distribution at population size $z$.
Because, by Assumption (A3), $\mbE\left(\xi(z)^{l}\right)$ is finite for all $l \geq 0$, we then have 
\[
\mbP(A_{2,n}) \geq 1 - \lceil C(i,q)\log n \rceil \,\mbP(X^\uparrow > n^{1/4} ) \geq 1 - \frac{C z^{r+\nu}\, \mbE\left(\xi(z)^{r+\nu}\right)}{\rho\,v_z n^{r/4-1}}.
\]
Selecting $r > 8$ we see that $\mbP(A_{2,n}) \to 1$. To show $\mbP(A_{3,n}) \to 1$ let $X^{(k)}$ be such that $\mbP(X^{(k)}=j) = Q_{zj} \frac{\bm{e}_j Q^{k-1} \bm{1}}{\bm{e}_z Q^k \bm{1}}$. 
Observe that for all $k< C(z,q)\log n$
\[
\mbP(X^{(k)}=j) \leq  \frac{Q_{zj}}{\bm{e}_z Q^{C(z,q)\log n} \bm{1}} \leq \frac{Q_{zj}}{v_z n^{-K(z,q)}},
\]
where the last equality uses \eqref{Qn} and therefore holds for $n$ sufficiently large, and $K(z,q)$ is a positive constant that depends on $z$ and $q$. Then, following the same arguments as before, we obtain
\[
\mbP(A_{3,n}) \geq 1- \frac{z^r\mbE\left(\xi(z)^{r}\right)}{v_z n^{r/4-K(z,q)-1}},
\]
and we can choose $r$ large enough to ensure $\mbP(A_{3,n}) \to 1$. \end{pro}


\begin{lem}\label{le2}
For any initial state $i\geq1$ and every $z\in\N$, conditional on $Z_0=i$, we have $\widehat Y_n^\uparrow(z) \stackrel{d}{\to} Z$.
\end{lem}

\begin{pro}
We fix $z$ and write $\widehat Y_n^\uparrow(z)=(z\,u_z v_z n)^{-1/2}{M_n},$ where
$$M_0:=0,\quad M_n:=\sum_{i=1}^n \{Z^\uparrow_{i}-z{m}^{\uparrow}(z)\} \ind{Z^\uparrow_{i-1}=z},\quad n\geq 1,$$forms a mean-0 martingale in $L^2$ with increments $Y_k:= \{Z^\uparrow_{k}-z{m}^{\uparrow}(z)\} \ind{Z^\uparrow_{k-1}=z},$ $1\leq k\leq n$.
We verify Conditions \eqref{clt-cond1} and \eqref{clt-cond2} to use the martingale central limit theorem (MCLT) in Section \ref{mclt}.
Let $Z(z)$ be a random variable with the same distribution as $Z^\uparrow_{k}$ conditional on $Z^\uparrow_{k-1}=z$, and observe that $\mbE\left(Z(z)\right)=z{m}^{\uparrow}(z)$ and $\mbV\left(Z(z)\right)=z^2\,{\sigma^2}^{\uparrow}(z)$.
We have, for any $i\geq 1$, $$\mbE_i(Y_k^2\,|\,\mathcal{F}_{k-1})=\mbE_i\left(\{Z^\uparrow_{k}-z{m}^{\uparrow}(z)\}^2 \ind{Z^\uparrow_{k-1}=z}\,|\,Z^\uparrow_{k-1}\right)=\mbV\left(Z(z)\right)\,\ind{Z^\uparrow_{k-1}=z\,|\,Z^\uparrow_{0}=i}.$$Therefore, from the ergodic theorem for the positive recurrent Markov chain $\{ Z^\uparrow_\ell\}$,
$$\lim_{n\to\infty}n^{-1} \sum_{k=1}^n \mbE_i(Y_k^2\,|\,\mathcal{F}_{k-1}) =z^2\,{\sigma^2}^{\uparrow}(z) u_zv_z$$in probability.  Condition \eqref{clt-cond1}  is therefore satisfied with $\Gamma=z^2\,{\sigma^2}^{\uparrow}(z) u_zv_z$. In addition, 
\begin{eqnarray*}\mbE_i(Y_k^2\,\ind{|Y_k|\geq \varepsilon \sqrt{n}}|\,\mathcal{F}_{k-1})&=&\mbE_i\left(\{Z^\uparrow_{k}-z{m}^{\uparrow}(z)\}^2 \ind{Z^\uparrow_{k-1}=z}\,\ind{|Z^\uparrow_{k}-z{m}^{\uparrow}(z)|\ind{Z^\uparrow_{k-1}=z}\geq \varepsilon \sqrt{n}}|\,Z^\uparrow_{k-1}\right)\\&=&\mbE\left(\{Z(z)-z{m}^{\uparrow}(z)\}^2 \,\ind{|Z(z)-z{m}^{\uparrow}(z)|\geq \varepsilon \sqrt{n}}\right)\,\ind{Z^\uparrow_{k-1}=z\,|\,Z^\uparrow_{0}=i}.
\end{eqnarray*}
Since $ \mbE\left(\{Z(z)-z{m}^{\uparrow}(z)\}^2 \,\ind{|Z(z)-z{m}^{\uparrow}(z)|\geq \varepsilon \sqrt{n}}\right)\rightarrow 0$ in probability as $n\to\infty$, Condition  \eqref{clt-cond2} is satisfied. It then follows from the MCLT that $n^{-1/2}\,M_n$ converges weakly to a centred normal random variable with variance $\Gamma$, and therefore $\widehat Y_n^\uparrow(z)$ converges weakly to a centred normal random variable with variance $\Gamma/(z u_z v_z)=z {\sigma^2}^{\uparrow}(z),$ which completes the proof. 

\end{pro}

We now have all the ingredients to complete the proof of Theorem \ref{thm:psdbp-Qcons and Qnorm}.

\medskip
\begin{Prf}[Theorem \ref{thm:psdbp-Qcons and Qnorm}.]
%
Applying Lemmas \ref{le1} and \ref{le2} and \cite[Theorem 25.4]{Billingsley}  then leads to $\widehat Y_n^{(n)}(z) \stackrel{d}{\to} Z$, which proves \eqref{mz_Qnorm}.
Since, conditional on $Z_n>0$, the variance of $\hat{m}_n(z)-{m}^{\uparrow}(z)$ vanishes asymptotically, \eqref{mz_Qcons} follows.

To show that ${m}^{\uparrow}(z)$ is finite for every $z$, we express it in terms of the original process as ${m}^{\uparrow}(z)=(z \,\rho\, v_z)^{-1}\sum_{k\geq 1} k v_k Q_{zk}.$ 
Assumptions \ref{cond:irreducible-PSDBP}--\ref{cond:bounded-moments-xi} and Lemma \ref{Goss2} therefore ensure that ${m}^{\uparrow}(z)<\infty$ since,
for any $\nu>\nu^*$, there exists $C$ such that  $v_k/k^\nu < C$ for all $k$, hence
 $$\sum_{k\geq 1} k v_k Q_{zk}=\sum_{k\geq 1} k^{\nu+1} \dfrac{v_k}{k^{\nu}} Q_{zk}\leq C \sum_{k\geq 1} k^{\nu+1}  Q_{zk}=C \,\mbE(Z_1^{\nu+1}\,|\,Z_0=z)<\infty.$$A similar argument shows that ${\sigma^2}^{\uparrow}(z)<\infty$.

Finally, let us fix $z_1\neq z_2$ and study the asymptotic properties of 
$n^{1/2}(\hat{m}_{n}(z_1),\hat{m}_{n}(z_2))$ conditional on $Z_n>0$. Like in the proof of Lemma \ref{le2}, it mainly suffices to focus on the asymptotic normality of $n^{-1/2}(M_{n}(z_1),M_{n}(z_2))$, where
$$(M_{n}(z_1),M_{n}(z_2)):=\left(\sum_{i=1}^n \{Z^\uparrow_{i}-z_1{m}^{\uparrow}(z_1)\} \ind{Z^\uparrow_{i-1}=z_1},\sum_{i=1}^n \{Z^\uparrow_{i}-z_2{m}^{\uparrow}(z_2)\} \ind{Z^\uparrow_{i-1}=z_2}\right),\quad n\geq 1,$$ is a centred bivariate martingale in $L^2$  with increments $$Y_k:= \left(\{Z^\uparrow_{k}-z_1{m}^{\uparrow}(z_1)\} \ind{Z^\uparrow_{k-1}=z_1},\{Z^\uparrow_{k}-z_2{m}^{\uparrow}(z_2)\} \ind{Z^\uparrow_{k-1}=z_2}\right), \quad k\geq 1.$$Since, for $z_1\neq z_2$,
$$\mbE_i(Y_{k,1}Y_{k,2}\,|\,\mathcal{F}_{k-1})=\mbE_i\left(\{Z^\uparrow_{k}-z_1{m}^{\uparrow}(z_1)\} \ind{Z^\uparrow_{k-1}=z_1}\{Z^\uparrow_{k}-z_2{m}^{\uparrow}(z_2)\}^2 \ind{Z^\uparrow_{k-1}=z_2}\,|\,Z^\uparrow_{k-1}\right)=0,$$ the matrix $\Gamma$ in \eqref{clt-cond1} is diagonal, which proves asymptotic uncorrelation. 
\end{Prf}

\begin{Prf}[Corollary \ref{cor:psdbp-Qcons and Qnorm}.]
The fact that $\hat{m}_n(z)$ satisfies Definition \ref{def:Q-} (i) follows from Lemmas \ref{proof1} and \ref{le2}, and the fact that it satisfies Definition \ref{def:Q-} (ii) is straightforward from Theorem \ref{thm:psdbp-Qcons and Qnorm}.
\end{Prf}

\subsection{Proofs of the results in Section \ref{sec:Ccons}}\label{pr_GW}

\begin{Prf}[Theorem \ref{thm:GW-Qcons and Qnorm}.]
Conditionally on $Z_n>0$,
$$|\hat{m}_n-1|=\left|\dfrac{Z_n-Z_0}{\sum_{i=1}^{n}Z_{i-1}}\right|\leq
\dfrac{Z_n+Z_0}{n}.
$$
%
Therefore, on the set $\{Z_n>0\}$, for any $q<1$ and $\varepsilon>0$,
$$\{\varepsilon<n^q\,|\hat{m}_n-1|\}\subseteq \left\{\varepsilon<\frac{Z_n+Z_0}{n^{1-q}}\right\}\subseteq \left\{\frac{\varepsilon}{2}<\frac{Z_n}{n^{1-q}}\right\}\cup \left\{\frac{\varepsilon}{2}<\frac{Z_0}{n^{1-q}}\right\},$$
and as a consequence,
$$\mbP_{\pi}(n^q\,|\hat{m}_n-1|>\varepsilon\,|\,Z_n>0)\leq \mbP_{\pi}\left(\frac{\varepsilon}{2}<\frac{Z_n}{n^{1-q}}\,\Big|\,Z_n>0\right)+\mbP_{\pi}\left(\frac{\varepsilon}{2}<\frac{Z_0}{n^{1-q}}\,\Big|\,Z_n>0\right).$$
As $n\rightarrow\infty$, the first term in the right-hand-side converges to 0 because, conditionally on $Z_n>0$, the  law of $Z_n$ converges to the quasi-stationary distribution $\vc u$ (whose existence is guaranteed under the assumptions of the theorem), and therefore the law of $Z_n/n^{1-q}$ converges to the distribution of a degenerate variable at 0. The convergence of the second term in the right-hand-side follows from similar arguments by considering that, conditionally on $Z_n>0$, the distribution of $Z_0$ converges to the size-biased distribution of $\vc\pi$. 
\end{Prf}

\begin{Prf}[Corollary \ref{cor:GW-Qcons}]
In the subcritical GW case, we can interpret the universal limit 1 of $\hat{m}_{n}$ in terms of the $Q$-process.  Indeed, in that case, $\rho=m$, and the eigenvector $\vc v$ takes the simple form $v_i=C i$ for a positive constant $C$, and $i\geq 1$. 
Let $\vc w=(w_i)$ be the size-biased distribution associated with the stationary distribution of the $Q$-process, that is, $w_i=i u_i v_i/(\sum_j j u_j v_j)= i^2 u_i/(\sum_j j^2 u_j)$. We then have
\begin{eqnarray*}
f(Q^\uparrow):=\sum_i w_i \dfrac{m^{\uparrow}(i)}{i}&=& (\sum_j j^2 u_j)^{-1}\sum_i  i^2 u_i \dfrac{1}{i}\sum_{k} Q^{\uparrow}_{ik}\, k\\&=&(\sum_j j^2 u_j)^{-1}\sum_i  i u_i \sum_{k} \dfrac{Q_{ik}}{m}\, \dfrac{k^2}{i}\\&=&(\sum_j j^2 u_j)^{-1}\sum_{k}\dfrac{m \,u_k}{m}k^2\\&=&1,
\end{eqnarray*}
where we used the fact that $\vc u Q=m \vc u$.   So the value 1 in a way represents the asymptotic mean offspring per individual in the $Q$-process.  Note that in addition, $m(i)=i \,m$, and therefore $f(Q)=\sum_i w_i ({m(i)}/{i})=m$. We have then established (i) of Definition \ref{def:Q-}.

By Theorem \ref{thm:GW-Qcons and Qnorm}, conditionally on survival, the classical estimator $\hat{m}_{n}$ therefore converges  to the equivalent quantity of the mean offspring $m$ in the $Q$-process, and is thus $Q$-consistent according to (ii) of Definition \ref{def:Q-}.

\end{Prf}


\medskip

In order to prove asymptotic properties for the estimators $\tilde{m}_n$ and $\bar{m}_n$ of the mean offspring $m=\mbE(\xi)$, it is useful to interpret the process $\{X_n:=Z^{\uparrow}_n -1\}$ as a GW process with immigration, where the offspring distribution is $\xi$ and the distribution of the number of immigrants is  \mbox{SB$(\xi)-1$}. Following the notation in \cite{heyde1972estimation,heyde1972estimation_corr}, the mean number of immigrants is then given by 
\begin{equation}\label{llam}\lambda:=\mbE(\text{SB}(\xi)-1)=\dfrac{\mbV(\xi)}{m}+m-1,\end{equation}
the stationary mean of $\{X_n\}$ is
\begin{equation}\label{mmu}\mu:=\dfrac{\lambda}{1-m}=\dfrac{\mbV(\xi)}{m(1-m)}-1,\end{equation}and we let \begin{equation}\label{cc} c^2:=\mbV(\text{SB}(\xi)-1)+\mbV(\xi)\,\mu=\dfrac{\mbE(\xi^3)}{m}-\dfrac{\mbV(\xi)^2(1-2m)}{m^2\,(1-m)}-3\mbV(\xi)-m^2.\end{equation}

We first turn our attention to the consistency and asymptotic normality of the estimator
\[
\tilde m_n := \frac{\sum_{i=1}^n (Z_i - b )}{\sum_{i=1}^n(Z_{i-1} + a)},
\]
where we assume that 
\[
\frac{\mbV(\xi)}{m} = am+b.
\]

We start by establishing a few preliminary technical results.
We say two probability mass functions $p_1(\cdot)$ and $p_2(\cdot)$ satisfy the monotone likelihood ratio property if, for every $k > j$ we have 
\begin{equation}\label{mlr}
\frac{p_1(k)}{p_2(k)} \geq \frac{p_1(j)}{p_2(j)}. 
\end{equation}
In addition, if \eqref{mlr} holds and $X_1 \sim p_1$ and $X_2 \sim p_2$, then $X_1$ stochastically dominates $X_2$. A Markov chain $\{ Z_n \}$ is called stochastically monotone if for any $k>\ell$, and any $n\geq 1$, $(Z_n\,|\,Z_{n-1}=k)$ stochastically dominates $(Z_n\,|\,Z_{n-1}=\ell)$.

\begin{lem}
Suppose $\{ Z_n \}$ is stochastically monotone, then for all $i>0$, $\ell < n$ and $j < k$ we have 
\begin{equation}\label{mlrp1}
\frac{\mbP_i(Z_1 = k | Z_n>0)}{\mbP_i(Z_1 = k | Z_\ell>0)} \geq \frac{\mbP_i(Z_1 = j | Z_n>0)}{\mbP_i(Z_1 = j | Z_\ell>0)}
\end{equation}
\end{lem}
\begin{pro} Observe that \eqref{mlrp1} is equivalent to 
\[
\frac{Q_{ik}\frac{ \bm{e}^\top_k Q^{n-1} \bm{1}}{\bm{e}^\top_i Q^{n} \bm{1}}}{Q_{ik}\frac{ \bm{e}^\top_k Q^{\ell-1} \bm{1}}{\bm{e}^\top_i Q^{\ell} \bm{1}}} \geq \frac{Q_{ij}\frac{ \bm{e}^\top_j Q^{n-1} \bm{1}}{\bm{e}^\top_i Q^{n} \bm{1}}}{Q_{ij}\frac{ \bm{e}^\top_j Q^{\ell-1} \bm{1}}{\bm{e}^\top_i Q^{\ell} \bm{1}}},
\]
which is equivalent to 
\[
\frac{ \bm{e}^\top_k Q^{n-1} \bm{1}}{ \bm{e}^\top_k Q^{\ell-1} \bm{1}} \geq \frac{ \bm{e}^\top_j Q^{n-1} \bm{1}}{ \bm{e}^\top_jQ^{\ell-1} \bm{1}},
\]
that is,
\begin{equation}\label{mlrp2}
\mbP_k(Z_{n-1} > 0 | Z_{\ell-1} >0) \geq \mbP_j(Z_{n-1} >0 | Z_{\ell-1}>0 ).
\end{equation}
However, due to stochastic monotonicity, it is easily seen that $(Z_{\ell-1} | Z_{\ell-1}>0, Z_0=k)$ stochastically dominates $(Z_{\ell-1} | Z_{\ell-1}
 >0, Z_0 =j)$ (for instance by verifying the monotone likelihood ratio property). Combining this with stochastic monotonicity then implies \eqref{mlrp2}.
 \end{pro}

\begin{cor}\label{cor_st_mon}
If $\{Z_n \}$ is stochastically monotone, then for all $i \geq 1$ and $\ell < n < n+k$ we have that $(Z_\ell | Z_n >0, Z_0= i)$ is stochastically dominated by $(Z_\ell | Z_{n+k} >0, Z_0 =i)$. 
In particular, $(Z_\ell |  Z_n >0, Z_0= i)$ is stochastically dominated by $(Z^\uparrow_\ell | Z_0= i)$.
\end{cor}

We place ourselves in the probability space of a MEXIT coupling $\{(\widehat{Z}_\ell^{(n)},\widehat{Z}^\uparrow_\ell )\}$ defined in Section \ref{sec:coupling}, and we simply write $\mbP$ for the associated measure $\widehat{\mbP}^{(n,\uparrow)}(\cdot)$.
Recall that $\tau_n$ is the uncoupling time defined in \eqref{uncoupling} and $C(i,q)$ is the constant defined in Theorem \ref{coupling1} (iii).
We define the events
\begin{align*}
A_{1,n} &:= \{ \tau_n > n - C(i,q) \log n \}, \\
A_{2,q,n} &:= \{ \hat{Z}^\uparrow_\ell \leq n^{q}\;\;\text{for all } \ell:   n - C(i,q) \log n \leq \ell\leq n \}, \quad q>0,\\ A_{3,q,n} &:= \{ \hat{Z}^{(n)}_\ell \leq n^{q}\;\;\text{for all } \ell:   n - C(i,q) \log n \leq \ell\leq n \}, \quad q>0.
\end{align*}

\begin{lem}\label{lemAAA}Under Assumptions (A1)--(A3),
for all $i\geq1$ and $q>0$, we have $\mbP_i(A_{1,n}) \to 1$, $\mbP_i(A_{2,q,n}) \to 1$, and $\mbP_i(A_{3,q,n}) \to 1$ as $n\to\infty$.

\end{lem}

\begin{pro}
From Theorem \ref{coupling1} (iii) we have $\mbP_i(A_{1,n}) \to 1$. 
To show $\mbP_i(A_{2,q,n}) \to 1$, we observe that, when $m<1$,
\[
\mbE_i(Z^\uparrow_\ell) = i m^\ell + \dfrac{\mbV(\xi)}{m}\sum_{j=1}^{\ell-1} m^{j} \leq  i + \dfrac{\mbV(\xi)}{m(1-m)},
\]
for all $\ell \geq 0$. Thus, applying Markov's inequality, we obtain 
\[
\mbP_i( Z_\ell^\uparrow > n^{q} ) \leq n^{-q} \left(i + \dfrac{\mbV(\xi)}{m(1-m)}\right),
\]
so that $$\mbP_i(A_{2,q,n}) \geq 1- C(i,q) (\log n)\,n^{-q} \left(i + \dfrac{\mbV(\xi)}{m(1-m)}\right) \to 1.$$ 
To show  $\mbP_i(A_{3,q,n}) \to 1$ we observe that Corollary \ref{cor_st_mon} implies $\mbE_i(Z^\uparrow_\ell) \leq \mbE_i(Z_\ell | Z_n >0)$. We then apply the same arguments.
\end{pro}

We now show an equivalent of Lemma \ref{proof1}. Let $m_\pi$ be the mean of the stationary distribution of the $Q$-process; it is given by $m_\pi=\mu+1$, where $\mu$ is given in \eqref{mmu}. Under Assumption (A4), we have $m_\pi={(am+b)}/{(1-m)}$.

 \begin{lem}\label{proof2}
 Under Assumptions (A1)--(A3),
 for all initial population size $i\geq 1$ and any $\varepsilon>0$, we have 
 \[
 \mbP_i \left( \left| \frac{\sum_{\ell=1}^n  Z_\ell }{n} - m_\pi \right| > \varepsilon 
 , \bigg| \, Z_n >0\right) \to 0
 \]
 as $n \to \infty$.
 \end{lem}
 
  \begin{pro}
 Observe that, on the paths where $A_{1,n}, A_{2,q,n}$, and $A_{3,q,n}$ happen for $0<q<1$, 
\begin{equation}\label{diff_zero2}
 n^{-1} \left| \sum_{\ell=1}^n  \hat{Z}_\ell^{(n)} - \sum_{\ell=1}^n  \hat{Z}_\ell^{\uparrow}  \right| \leq \frac{C(i,q)(\log n)\,n^q }{n}.
\end{equation}
Together with Lemma \ref{lemAAA}, this implies that the left-hand-side of \eqref{diff_zero2} converges to $0$ in probability as $n \to \infty$. 
From the ergodic theorem for the Markov chain $\{ Z^\uparrow_\ell\}$,
we have 
$n^{-1} \sum_{\ell=1}^n  Z_\ell^{\uparrow} \to m_\pi$ in probability. The result then follows from \cite[Theorem 25.4]{Billingsley}.
 \end{pro}
 
To establish consistency and asymptotic normality of $\tilde m_n$ in the subcritical case  $m<1$, we proceed in a similar way to Section \ref{pr_PSDBP}.  
For $n\geq1$, we define the random variables
$$\widehat W_n^{(n)}:=n^{-1/2} \,\left(\sum_{i=1}^n Z_i^{(n)}-b-m\,Z_{i-1}^{(n)}-m\,a\right),\; \widehat W_n^\uparrow=n^{-1/2} \,\left(\sum_{i=1}^n Z_i^\uparrow-b-m\,Z_{i-1}^\uparrow-m\,a\right).$$ 
 
 \medskip
 \begin{lem}\label{le1b}
Under Assumptions (A1)--(A3), for any initial population size $i\geq1$, conditional on $Z_0=i$, we have $\widehat W_n^{(n)} - \widehat W_n^\uparrow \stackrel{P}{\to} 0$.
\end{lem} 

\begin{pro}The proof follows the same arguments as that of Lemma \ref{proof2} but now with $0<q<1/2$.
\end{pro}

\medskip

Under Assumption (A4), the constant $c^2$ in \eqref{cc} takes the  more explicit form
\begin{equation}
\label{cc2}c^2=\dfrac{\mbE(\xi^3)}{m}+\dfrac{m(am+b)^2}{1-m}-m(am+b)-[(a+1)m+b]^2.
\end{equation}

\begin{lem}\label{le2b}
Under Assumptions (A1)--(A3), for any initial population size $i\geq1$, conditional on $Z_0=i$, we have $\widehat W_n^\uparrow \stackrel{d}{\to} Z$,
where $Z\sim N(0,c^2)$ with $c^2$ given in \eqref{cc2}.
\end{lem}

\begin{pro}
The proof follows the same arguments as those used in the proof of  \cite[Theorem 3]{heyde1972estimation} (with further corrections in \cite{heyde1972estimation_corr}). The idea is to write  $$Y_i:= Z^\uparrow_{i}-\mbE(Z^\uparrow_{i}\,|\,\mathcal{F}_{i-1})=Z^\uparrow_{i}-b-Z^\uparrow_{i-1}{m}-ma,$$and $M_n:=\sum_{i=1}^n Y_i$, and note that
$\{(M_n,\mathcal{F}_n)\}_{n\geq 1}$ is a mean-0 martingale. We can then show that a slight variant of Condition \eqref{clt-cond1} is satisfied with $\Gamma=c^2$ (see \cite[pp. 246--247]{heyde1972estimation}) and a slight variant of Condition \eqref{clt-cond2} is also satisfied (see \cite[pp. 573]{heyde1972estimation_corr}). By the MCLT, we then have
$$\widehat W_n^\uparrow=n^{-1/2}\,M_n \stackrel{d}{\to} Z,$$
where $Z\sim N(0,c^2)$.

\end{pro}

\begin{Prf}[Proposition \ref{Ccons1}]
In the subcritical case $m<1$, conditionally on $Z_n>0$, we have
$$\sqrt{n}(\tilde m_n-m)=\widehat W_n^{(n)} \,\frac{n}{\sum_{\ell=1}^{n}(Z^{(n)}_{\ell-1}+a)}.$$
The  asymptotic normality result in \eqref{m_Cnorm1} then follows from Lemma \ref{proof2} with $m_\pi+a=(am+b)/(1-m)+a=(a+b)/(1-m)>0$, Lemmas  \ref{le1b} and \ref{le2b},  \cite[Theorem 25.4]{Billingsley}, and Slutsky's theorem, with
\begin{equation}\label{nut2}\tilde{\nu}^{2}:=c^2 \dfrac{(1-m)^2}{(a+b)^2},\end{equation}where $c^2$ is given in \eqref{cc2}.
Since the variance of $\tilde m_n-m$ vanishes asymptotically, $\tilde m_n$ is $C$-consistent for $m$.

In the supercritical case $m>1$, by Kesten-Stigum theorem, $Z_n/m^n$ converges a.s. to a random variable $W$ which is positive on the paths of survival. As a consequence, using Toeplitz lemma,
we have $m^{-n}\sum_{\ell=0}^n Z_\ell\to m W/(m-1)$ a.s., and conditional on $Z_0=i$ and on survival of the process,
\begin{eqnarray*}\frac{\sum_{\ell=1}^n (Z_\ell - b )}{\sum_{\ell=1}^n(Z_{\ell-1} + a)}&=&\dfrac{(\sum_{\ell=0}^n Z_\ell) -(i+n b)}{(\sum_{\ell=0}^{n-1} Z_\ell)+n a}\\&=&m\dfrac{m^{-n}[(\sum_{\ell=0}^n Z_\ell) -(i+n b)]}{m^{-(n-1)}[(\sum_{\ell=0}^{n-1} Z_\ell)+n a]}\\&\to& m\dfrac{m W/(m-1)}{m W/(m-1)}=m,\end{eqnarray*}which implies $C$-consistency of $\tilde{m}_n$ in the supercritical case.

In the critical case $m=1$, Yaglom's universal limit law states that, conditional on $Z_n>0$, $Z_n/n$ converges weakly to an exponential random variable $X$. We obtain $C$-consistency of $\tilde{m}_n$ in the critical case using the same argument as above by noting that $n^{-2}\sum_{\ell=0}^n Z_\ell\stackrel{d}{\to} X/2$.

\end{Prf}

\bigskip

Finally, we study the asymptotic properties of the estimators
\[
\bar m_n :=  1-\dfrac{1}{2}\dfrac{\sum_{i=1}^{n} (Z_i-Z_{i-1})^2}{\sum_{i=1}^{n}(Z_{i-1}-\bar{Z}_n)^2}
\]of the mean offspring $m$, and
$\bar{\sigma}^2_n:=\bar{m}_n(1-\bar{m}_n)\bar{Z}_n$ of the offspring variance $\sigma^2$,
where $\bar{Z}_n:=n^{-1}\,{\sum_{i=1}^n Z_{i-1}}$. We note from \cite{heyde1972estimation} that $\bar m_n$ exhibits the same asymptotic properties as the statistic \begin{equation}\label{mp}\bar m_n'=\dfrac{\sum_{i=1}^{n}(Z_i-\bar{Z}_n)(Z_{i-1}-\bar{Z}_n)}{\sum_{i=1}^{n}(Z_{i-1}-\bar{Z}_n)^2},\end{equation}so we prove our results for $\bar m_n'$.
We start by studying the asymptotic properties of the denominator  in \eqref{mp}.

 \begin{lem}\label{proof3}
 Under Assumptions (A1)--(A4),
 for all $i\geq 1$ and $\varepsilon>0$, we have 
 \[
 \mbP_i \left( \left| \frac{\sum_{\ell=1}^{n}(Z_{\ell-1}-\bar{Z}_n)^2}{n} - (1-m^2)^{-1} c^2 \right| > \varepsilon 
 , \Big| \, Z_n >0\right) \to 0
 \]
 as $n \to \infty$, where $c^2$ is given by \eqref{cc}.
 \end{lem}
 
 \begin{pro}
  Observe that, on the paths where $A_{1,n}, A_{2,q,n}$, and $A_{3,q,n}$ happen for $0<q<1/2$, 
\begin{eqnarray*}
\lefteqn{n^{-1} \left| \sum_{\ell=1}^n  \left\{(Z_{\ell-1}^{(n)}-\bar{Z}_n^{(n)})^2 -  (Z_{\ell-1}^{\uparrow}-\bar{Z}_n^{\uparrow})^2\right\} \right|}\\&=&\left|(\bar{Z}_n^{(n)})^2-(\bar{Z}_n^{\uparrow})^2+ n^{-1} \sum_{\ell=1}^n  \{(Z_{\ell-1}^{(n)})^2-(Z_{\ell-1}^{\uparrow})^2\} +2 \bar{Z}_n^{\uparrow}  \,n^{-1}\sum_{\ell=1}^nZ_{\ell-1}^{\uparrow} - 2 \bar{Z}_n^{(n)} \, n^{-1}\sum_{\ell=1}^n Z_{\ell-1}^{(n)}
 \right| \\ &=&\left|(\bar{Z}_n^{\uparrow})^2-(\bar{Z}_n^{(n)})^2+ n^{-1} \sum_{\ell=1}^n  \{(Z_{\ell-1}^{(n)})^2-(Z_{\ell-1}^{\uparrow})^2\} \right|\\&\leq& \left|(\bar{Z}_n^{\uparrow})^2-(\bar{Z}_n^{(n)})^2\right| +  \frac{C(i,q)(\log n)\,n^{2q} }{n};
\end{eqnarray*}(note that to keep the presentation light, here we omitted the \,$\hat{}$\, symbol on the coupled random variables).
By Lemma \ref{proof2}, the ergodic theorem for the Markov chain $\{ Z^\uparrow_\ell\}$, and the continuous mapping theorem, we have $\left|(\bar{Z}_n^{\uparrow})^2-(\bar{Z}_n^{(n)})^2\right|\rightarrow 0$ in probability.  By Lemma \ref{lemAAA}, the left-hand-side of the above equation therefore converges to $0$ in probability as $n \to \infty$.
 In addition, by \cite[Equation (3.13)]{heyde1972estimation}, we have $n^{-1}\sum_{\ell=1}^{n}(Z_{\ell-1}^{\uparrow}-\bar{Z}_n^{\uparrow})^2\rightarrow  (1-m^2)^{-1} c^2$ in probability.  The result then follows from  Lemma \ref{lemAAA} and
\cite[Theorem 25.4]{Billingsley}.
 
 \end{pro}
 
 To deal with the numerator in $\bar{m}'_n-m$, we define the random variables
\begin{eqnarray*}\widehat V_n^{(n)}&:=&n^{-1/2} \,\left(\sum_{i=1}^n (Z_{i-1}^{(n)}-\bar{Z}_n^{(n)})\{(Z_i^{(n)}-\bar{Z}_n^{(n)})-m\,(Z_{i-1}^{(n)}-\bar{Z}_n^{(n)})\}\right),\\ V_n^\uparrow&:=&n^{-1/2} \,\left(\sum_{i=1}^n (Z_{i-1}^\uparrow-\bar{Z}_n^\uparrow)\{(Z_i^\uparrow-\bar{Z}_n^\uparrow)-m\,(Z_{i-1}^\uparrow-\bar{Z}_n^\uparrow)\}\right).\end{eqnarray*}

 \begin{lem}\label{le1c}Under the assumptions of Proposition \ref{Ccons2},
for any initial state $i\geq1$, conditional on $Z_0=i$, we have $\widehat V_n^{(n)} - \widehat V_n^\uparrow \stackrel{P}{\to} 0$.
\end{lem} 

\begin{pro}
The proof follows the same arguments as that of Lemma \ref{proof3} but now with $0<q<1/4$.

\end{pro}

We introduce the constant
\begin{equation} \label{B}B^2:=c^4(1-m^2)^{-1}+\sigma^2(1-m^3)^{-1}\{\mbE[(\text{SB}(\xi)-1-\lambda)^3]+\mu\mbE[(\xi-m)^3]+3m\,\sigma^2 c^2(1-m^2)^{-1}\},
\end{equation}
with $\lambda, \mu,$ and  $c^2$  given in \eqref{llam}, \eqref{mmu}, and \eqref{cc}, respectively (the constant $B^2$ is obtained from  \cite[Theorem B and Theorem 3]{heyde1972estimation}, with further corrections in \cite{heyde1972estimation_corr}).

\begin{lem}\label{le2c}Under the assumptions of Proposition \ref{Ccons2},
for any initial state $i\geq1$, conditional on $Z_0=i$, we have $\widehat V_n^\uparrow \stackrel{d}{\to} Z$,
where $Z\sim N(0,B^2)$ with $B^2$ 
 given by \eqref{B}.
\end{lem}

\begin{pro}The proof follows from \cite[Theorem B (see also Theorem 2, and Theorem 3)]{heyde1972estimation} by interpreting the process $\{Z^{\uparrow}_n -1\}$ as a GW process with immigration (and offspring law $\xi$), where the law of the number of immigrants is  SB$(\xi)-1$.
\end{pro}

\begin{Prf}[Proposition \ref{Ccons2}]
Conditionally on $Z_n>0$, we have
$$\sqrt{n}(\bar m_n'-m)=\widehat V_n^{(n)} \,\frac{n}{\sum_{\ell=1}^{n}(Z^{(n)}_{\ell-1}-\bar{Z}_n^{(n)})^2}.$$
The asymptotic normality results in \eqref{m_Cnorm2}  
then follows from Lemmas \ref{proof3}, \ref{le1c}, \ref{le2c},  \cite[Theorem 25.4]{Billingsley}, and Slutsky's theorem, with
\begin{equation}\label{nu}
\bar{\nu}^{2}:=B^2(1-m^2)^2c^{-4},
\end{equation}
where $B^2$ is given in \eqref{B}.


Since, conditional on $Z_n>0$, the variance of $\bar m_n-m$ vanishes asymptotically, $\bar m_n$ is $C$-consistent for $m$. By \eqref{mmu}, $\mbV(\xi)=m(1-m)(\mu+1),$ where $\mu$ is the mean of the stationary distribution of $\{Z_n^\uparrow-1\}$. 
Therefore, by the fact that $\bar m_n$ is $C$-consistent for $m$, by Lemma \ref{proof2}, and by the continuous mapping theorem, we obtain that $\bar{\sigma}^2_n:=\bar{m}_n(1-\bar{m}_n)\bar{Z}_n$ is $C$-consistent for $\sigma^2$.

\end{Prf}
\subsection{Proofs of the results in Section \ref{sec:coupling}}

\subsubsection{Linear operator theory}\label{lot}

The proof of Theorem \ref{coupling1} relies on linear operator theory. Let $E$ be a Banach space. We say that $\lambda$ is an \emph{eigenvalue} of $Q$ on $E$ if there exists $\bm{v} \in E$ such that $Q \bm{v} = \lambda \bm{v}$, and we let $r_E(Q)$ denote the spectral radius of $Q$, that is, the supremum of $|\lambda|$ over all the eigenvalues $\lambda$.
The idea  is to show that the infinite matrix $Q$ displays properties more commonly associated with finite matrices. This essentially boils down to demonstrating that $Q$ is \emph{quasi-compact} on a carefully chosen Banach space. 
In short, $Q$ is quasi-compact on $E$ if there exists a decomposition 
\[
Q = \begin{bmatrix} 
Q_{1,1} & Q_{1,2} \\
Q_{2,1} & Q_{2,2}
\end{bmatrix},
\]
such that $Q_{1,1}$ is finite and $r_E(Q_{1,1}) > r_E(Q_{2,2})$ (see \cite{sasser} for a formal definition). Roughly speaking, this means the asymptotic behaviour of $Q^n$ is controlled by the finite matrix $Q_{1,1}$. 
In particular, for quasi-compact operators, there is a spectral gap between the largest and second largest eigenvalues on $E$ (see Lemma \ref{Goss2} below).

We let $ \lVert \bm{x} \rVert_{\infty,t} := \sup_{j \geq 1} t(j)^{-1} | x_j |$ for some function $t:\mathbb{N}\rightarrow\mathbb{R}$ and  $\lVert \bm{x} \rVert_{1,t}:= \sum_{j \geq 1} t(j) | x_j |$, and define the corresponding Banach spaces $l_t^\infty=\{ \bm{x} \in \mbR^{\mbN} : \lVert \bm{x} \rVert_{\infty,t} < \infty\}$ and $l^1_t=\{ \bm{x} \in \mbR^{\mbN} : \lVert \bm{x} \rVert_{1,t} < \infty\}$.
Define the \emph{matrix norm} of $Q$ in the Banach space $E$ as
\[
\lVert Q \rVert_E = \sup \{ \lVert Q\bm{v} \rVert_E : \lVert \bm{v} \rVert_E \leq 1 \}.
\]
We point out that $\lVert Q \rVert_{\infty,t} = \lVert Q^\top \rVert_{1,t}$, and that $l_{1,t}$ is referred to as the \emph{transpose} of $l_{\infty,t}$. 

Now suppose $t=t_0$ where $t_0$ is the identity operator (i.e. $t_0(x)=x$), and let $Q$ be the transition matrix of a subcritical GW branching process.
In this case it is well established (see for instance \cite[Chapter 3]{asmussen}) that $r_{1,t}(Q)=1 \neq \rho=m$ and that there exists a continuum of real-valued eigenvalues of $Q^\top$ on $l_{1,t}$ that lie between $m$ and $1$, whose corresponding eigenvectors are the quasi-stationary distributions of the process. 
To establish quasi-compactness for GW processes, and more generally for PSDBPs, we then need to select $t$ in such a way that, in the corresponding Banach space $E$, the continuum between $\rho$ and $r_{1,E}(Q)$ disappears (i.e so that $r_{1,E}(Q)=\rho=m$ in the GW case). 
Thankfully, for PSDBPs, this problem has already been addressed in \cite[Section 4]{Gosselin-2001}. The next lemma follows from \cite[Theorem 4.1 and Proposition 5.3]{Gosselin-2001}. We let $t_0^\nu$ be the function such that $t_0^\nu(x)= x^\nu$. 

\begin{lem}\label{Goss1}
Under assumptions (A1)--(A3), there exists $\nu^* \in \mbN$ such that for all $\nu > \nu^*$, $Q$ is quasi-compact on $l^\infty_{t_0^\nu}$.
\end{lem}

We are then able to use known properties of quasi-compact operators. In particular, the next lemma corresponds to \cite[Lemma 6.1]{Gosselin-2001}.
\begin{lem}\label{Goss2}
Under assumptions (A1)--(A3),  there exist $\nu^* \in \mbN$ and an operator $S$  such that
\[
Q^n = r(Q)^n \bm{v} \bm{u}^\top + S^n,
\]
where, for any $t=t_0^\nu$ with $\nu > \nu^*$, $r_{\infty,t}(S) < r_{\infty,t}(Q)$, and $\bm{u}$ and $\bm{v}$ are non-negative vectors such that $\bm{u} \in l^{1}_{t}$, $\bm{v} \in l_t^\infty$, and $\sum_{i=1}^\infty u_i v_i =1$.
\end{lem}
In addition, we have $( \lVert S^n \rVert_{\infty,t} )^{1/n} \to r_{\infty,t}(S)$ (Gelfand's theorem), which along with Lemma~\ref{Goss2} and the fact that $S$ is a bounded operator on $l^\infty_t$ (that is, $\lVert S \rVert_{\infty,t}<\infty$) implies that there exists $K<\infty$ and $\varepsilon >0$ such that, for any $n\geq1$,
\begin{equation}\label{Gelf}
\rho^{-n} \lVert S^n \rVert_{\infty,t} \leq K(1-\varepsilon)^n,
\end{equation}
with the same holding for $S^\top$ on the transpose space $l_{1,t}$.

\subsubsection{MEXIT coupling}\label{MEXIT}

For two probability measures $p_1(\cdot)$ and $p_2(\cdot)$ on $\mbN^\ell$, it is well known (see for instance \cite{thorisson1995coupling}) that there exists a coupling $(\hat X_1, \hat X_2)$ with marginal distributions $p_1(\cdot)$ and $p_2(\cdot)$, respectively, such that 
\[
\widehat{\mbP}(\widehat X_1 \neq \widehat X_2) = \frac{1}{2} \sum_{\bm{s} \in \mbN^\ell} | p_1(\bm{s}) - p_2(\bm{s}) |.
\]
This is the \emph{maximal coupling} that maximises $\widehat{\mbP}(\widehat X_1 = \widehat X_2)$. Now suppose $\{ \widehat X_{1,n} \}_{n \geq 0}$ and $\{ \widehat X_{2,n} \}_{n \geq 0}$ are Markov chains. Then the authors of \cite{ernst2019mexit} expand on this idea to demonstrate the existence of a \emph{MEXIT} coupling for Markov chains. The MEXIT coupling maximises the random variable 
\[
\tau := \min \{ \ell : \widehat X_{1,\ell} \neq \widehat X_{2, \ell} \},
\]
that is, $\tau$ under the MEXIT coupling stochastically dominates $\tau$ under any other coupling.
In particular, they demonstrate that, in the MEXIT coupling, for any (common) initial  state $i\geq1$ and all $\ell\geq 1$,
\begin{equation}\label{Mexc}
\widehat{\mbP}_i(\tau \leq \ell) = \frac{1}{2}\sum_{\bm{s} \in \mbN^\ell} | p_1^{(\ell)}(\bm{s}) - p_2^{(\ell)}(\bm{s}) |, 
\end{equation}
where $p_1^{(\ell)}(\bm{s})= \mbP_i( X_{1,u} = s_u, \, 1 \leq u \leq \ell)$ and $p_2^{(\ell)}(\bm{s})= \mbP_i( X_{2,u} = s_u, \, 1 \leq u \leq \ell)$. 

\subsubsection{Proof of Theorem \ref{coupling1}}

Let $\widehat{\mbP}_i^{(n,\uparrow)}(\cdot)$ denote a MEXIT coupling of $\{ Z_\ell ^{(n)}| Z_0=i \}_{1 \leq \ell \leq n}$ and $\{ Z^\uparrow_\ell | Z^\uparrow_0 = i \}_{1\leq \ell \leq n}.$
By \eqref{Mexc} we have, for every fixed $\ell\leq n$,
\begin{equation}\label{ce1}
\widehat{\mbP}_i^{(n,\uparrow)}( \tau_n \leq \ell) =\frac{1}{2}  \sum_{\bm{x} \in \mbN^{\ell}} |p_i^{(\ell,n)}(\bm{x})-p^{(\ell,\uparrow)}_i(\bm{x})|,
\end{equation}
where $p_i^{(\ell,n)}(\bm{x}) = \mbP_i( Z_j = x_j,  \, 1\leq j \leq \ell\, |\, Z_n >0)$ and $p_i^{(\ell,\uparrow)}(\bm{x}) = \mbP_i( Z^\uparrow_j = x_j, \,  1\leq j \leq \ell)$.
To bound the right-hand-side of \eqref{ce1}, we write
\begin{align}\nonumber
 &\sum_{\bm{x} \in \mbN^{\ell}} |p_i^{(\ell,n)}(\bm{x})-p^{(\ell,\uparrow)}_i(\bm{x})|  \\
&=  \sum_{\bm{x} \in \mathbb{N}^\ell } \left| \frac{Q_{ix_1} Q_{x_1x_2}\dots Q_{x_{\ell-1}x_\ell} \bm{e}^\top_{x_\ell} Q^{n-\ell} \bm{1}}{\bm{e}^\top_i Q^n \bm{1} } - \frac{Q_{ix_1}Q_{x_1x_2} \dots Q_{x_{\ell-1}{x_\ell}}v_{x_\ell}}{\rho^\ell v_i} \right|  \nonumber\\
&= \sum_{j =1}^\infty \bm{e}^\top_i Q^\ell \bm{e}_j \left| \frac{\bm{e}^\top_j Q^{n-\ell} \bm{1}}{\bm{e}^\top_i Q^n \bm{1}} - \frac{v_j}{\rho^\ell v_i} \right| \label{ce0}\\
&= \sum_{j=1}^\infty ( \rho^\ell u_j v_i + \bm{e}^\top_i S^\ell \bm{e}_j) \left| \frac{v_j \rho^{n-\ell} + \bm{e}^\top_j S^{n-\ell} \bm{1}}{v_i \rho^n + \bm{e}^\top_i S^n \bm{1} } - \frac{v_j}{\rho^\ell v_i}  \right| \nonumber\\
&=\ \sum_{j=1}^\infty ( \rho^\ell u_j v_i + \bm{e}^\top_i S^\ell \bm{e}_j) \left| \frac{ \rho^\ell v_i \bm{e}^\top_j S^{n-\ell} \bm{1}- v_j\bm{e}^\top_i S^n \bm{1} }{\rho^\ell v_i(v_i \rho^n + \bm{e}^\top_i S^n \bm{1} )}  \right|. \label{ce2}
\end{align}
Now, using H\"older's inequality in the first step and Equation \eqref{Gelf} in the last, for any $k \geq 1$, we have that there exists $K<\infty$ and $\varepsilon>0$ such that, for any function $t$ identified in Lemma~\ref{Goss2},
\begin{equation}\label{ce2a}
|\bm{e}^\top_i S^k \bm{e}_j | \leq \lVert (S^k)^\top \bm{e}_i \rVert_{1,t}\lVert \bm{e}_j \rVert_{\infty,t} \leq \lVert (S^k)^\top \rVert_{1,t} \lVert \bm{e}_i \rVert_{1,t}  \lVert \bm{e}_j \rVert_{\infty,t} \leq \frac{t(i)}{t(j)} K(1-\varepsilon)^k \rho^k,
\end{equation}
and, similarly,
\begin{equation}\label{ce2b}
|\bm{e}^\top_i S^k \bm{1} | \leq t(i) K(1-\varepsilon)^k\rho^k.
\end{equation}
For any $\eta > 0$ we can choose $N(i)$ large enough so that 
\[
t(i) K(1-\varepsilon)^n < \eta v_i
\]
for all $n > N(i)$. Now let $k=n-\ell$ (fixed). For $n>N(i)+k$ we have the upper bound 
\begin{align}
\text{\eqref{ce2}} &\leq 
 \sum_{j=1}^\infty ( \rho^\ell u_j v_i + \eta v_i \rho^\ell/t(j)) \left( \frac{ \rho^\ell v_i |\bm{e}^\top_j S^{k} \bm{1}| + v_j v_i \rho^n \eta }{\rho^\ell v_i(v_i \rho^n - v_i \rho^n \eta )} \right) \nonumber \\
&= \sum_{j=1}^\infty (u_j  + \eta/t(j)) \left( \frac{\rho^{-k}|\bm{e}_j S^k \bm{1}| + v_j \eta }{1 - \eta}  \right), \label{ce3u}
\end{align}
and the lower bound 
\begin{equation}
\text{\eqref{ce2}} \geq  \sum_{j=1}^\infty (u_j  - \eta/t(j)) \left( \frac{\rho^{-k}|\bm{e}_j S^k \bm{1}| - v_j \eta }{1 + \eta}  \right),\label{ce3l}
\end{equation}where we used the fact that, if $|B|\leq B^*$, then $|A|-B^*\leq |A|-|B|\leq |A-B|$.

Using the fact that $\sum_{j=1}^\infty u_j v_j =1$ we have 
\begin{equation}\label{ce4}
\text{\eqref{ce3u}} = \frac{\rho^{-k} }{1-\eta}\sum_{j=1}^\infty u_j\,|\bm{e}_j S^k \bm{1}| + \frac{\eta}{1-\eta} + \frac{\eta \rho^{-k}}{1-\eta} \sum_{j=1}^\infty \frac{|\bm{e}_j^\top S^{k} \bm{1}|}{t(j)} + \frac{\eta^2}{1-\eta}\sum_{j=1}^\infty \frac{v_j}{t(j)},
\end{equation}
with a similar equality for \eqref{ce3l}. 
If we can demonstrate that the sums in the third and fourth terms of \eqref{ce4} are finite, then the result follows by taking $\eta$ arbitrarily small. We consider the third term first. Using the fact that, for all $k \geq 1$,  
\begin{equation}\label{ce3b}
S^k = Q^k - \rho^k(\bm{v} \bm{u}^\top), \quad | \bm{e}_j^\top Q^k \bm{1} | \leq {1}, \quad \text{ and } \quad \bm{e}_j (\bm{v}\bm{u}^\top) \bm{1} = v_j,
\end{equation}
we obtain 
\begin{align}
\sum_{j=1}^\infty \left| \frac{\bm{e}_j^\top S^k \bm{1}}{t(j)} \right| &= \sum_{j=1}^\infty \left| \frac{\bm{e}_j^\top (Q^k - \rho^{k}(\bm{v} \bm{u}^\top)) \bm{1}}{t(j)} \right|  \leq  \sum_{j=1}^\infty \frac{1+ \rho^k v_j}{t(j)}. \label{ce5}
\end{align}
By Lemma \ref{Goss2} we have 
\[
\sup_{j \in \mbN} \frac{v_j}{j^{\nu^*}}<\infty,
\]
where $\nu^*$ is that of Lemma \ref{Goss2}. Because Lemma \ref{Goss2} holds for all $t = t_0^\nu$ with $\nu > \nu^*$, we can choose $\nu=\nu^*+2$, in which case \eqref{ce5} is finite. The same argument also applies to the fourth term in \eqref{ce4}. Taking $\eta$ arbitrarily small we thus obtain 
\[
\widehat{\mbP}_i^{(n,\uparrow)}( \tau_n \leq n-k) \to \frac{ \rho^{-k}}{2} \sum_{j=1}^\infty u_j |\bm{e}_j^\top S^k \bm{1}|
\]
as $n \to \infty$; this completes the  proof of \emph{(i)}.

Since $\widehat{\mbP}_i^{(n,\uparrow)}( \tau_n < \infty)=\widehat{\mbP}_i^{(n,\uparrow)}( \tau_n \leq n)$,  to establish \emph{(ii)}, in Equation \eqref{ce0}, we set $\ell:=n$, so  the term $\bm{e}_j^\top Q^{n-\ell} \bm{1}$ becomes 1. The same arguments can then be used to establish the result.

To establish \emph{(iii)}, we set $n-\ell := C(i,q) \log n$ ($=k$), and observe that, for any $q>0$, $\gamma>0$ and $t=t_0^\nu$, it is possible to select $C(i,q)$ $(=C(i,q,\nu,\gamma))$ large enough so that 
\begin{equation}\label{ce6}
t(i) K (1- \varepsilon)^{\lfloor C(i,q) \log n \rfloor} \leq \frac{\gamma v_i}{n^q},
\end{equation}
for all sufficiently large $n$. As a consequence, for sufficiently large $n$, we also have
\begin{equation}\label{ce6b}
t(i) K (1- \varepsilon)^n \leq \frac{\gamma v_i}{n^q}, \quad\textrm{and}\quad t(i) K (1- \varepsilon)^\ell \leq \frac{\gamma v_i}{n^q}.
\end{equation}
Combining \eqref{ce2}, \eqref{ce2a}, \eqref{ce2b}, \eqref{ce6}, and \eqref{ce6b}, we obtain
\begin{align}\nonumber
\widehat{\mbP}_i^{(n,\uparrow)}(\tau_n \leq \ell) &\leq \frac{1}{2}
\sum_{j=1}^\infty ( \rho^\ell u_j v_i + \frac{\gamma n^{-q} v_i \rho^\ell}{t(j)}) \left( \frac{ \rho^n v_i v_j \gamma n^{-q}+ v_j v_i \rho^n \gamma n^{-q} }{\rho^\ell v_i(v_i \rho^n - v_i \rho^n \gamma n^{-q} )} \right)
\\\label{2s}&\leq \frac{\gamma n^{-q}}{(1 -  \gamma n^{-q} )}\left(
\sum_{j=1}^\infty  u_j v_j  + \gamma n^{-q} \sum_{j=1}^\infty \frac{v_j}{t(j)}\right)
\\\nonumber
&\leq \frac{K_2 \gamma }{n^{q}},
\end{align}
for some constant $K_2$ and sufficiently large $n$.  Indeed, $\sum_{j=1}^\infty u_j v_j =1,$ and by choosing $t=t_0^{\nu^*+2}$, the second sum in \eqref{2s} is finite by the same argument as the one used to prove finiteness of \eqref{ce5}. The result then follows by choosing $\gamma$ sufficiently small so that $K_2 \gamma\leq 1$.

 \qed
 
 \appendix
 
 \section{Multilevel splitting method for the simulation of subcritical GW trajectories} 
 
 The numerical analysis of the estimators for the mean offspring of subcritical GW processes requires simulating long non-extinct trajectories of these processes. Due to the rapid extinction of subcritical GW processes, this is similar to rare event simulation, and can be done by adapting the multilevel splitting method (see for instance \cite{glasserman1999multilevel}) to our setting. 
 
 More precisely, to obtain a non-extinct trajectory of length $n$, we decompose the $n$ generations into sub-intervals containing $s$ generations (where $s$ has to be chosen optimally, and where the last sub-interval may be shorter than $s$). We start by simulating a first trajectory from generation $0$ until generation $s$. If that trajectory is extinct by generation $s$, we start again a new trajectory from time 0, and we repeat this step until we obtain a non-extinct trajectory at generation $s$. We then duplicate the trajectory at generation $s$ and simulate the two copies independently from generation $s$ to generation $2s$. If by generation $2s$ all trajectories are extinct, we start the whole process again from 0; otherwise, we keep duplicating each non-extinct trajectory (and leave the extinct ones) until the end of the next sub-interval, and so on, until we obtain at least one non-extinct trajectory at generation $n$; we then pick one of those at random as our sample trajectory. 
 
 The  simulated process can be seen as the output of a  ``macro GW process'' with binomial offspring distribution with mean $M(s)=2\,\mbP(Z_s>0)$, and the simulation is optimal when the macro process is critical; in practice we therefore choose $s$ to be the largest integer such that $M(s)\geq1$.

\section*{Acknowledgements}
Peter Braunsteins has conducted part of the work while supported by the Australian Research Council (ARC) Laureate Fellowship FL130100039. Sophie Hautphenne would like to thank the Australian Research Council (ARC) for support through her Discovery Early Career Researcher Award DE150101044. 
Carmen Minuesa's research has been supported by the Ministerio de Econom\'ia y Competitividad (grant MTM2015-70522-P), the Ministerio de Ciencia e Innovaci\'on (grant PID2019-108211GB-I00), the Junta de Extremadura (grants IB16099 and GR18103) and the Fondo Europeo de Desarrollo Regional. This research was initiated while Carmen Minuesa was a visiting postdoctoral researcher at The University of Melbourne, and she is grateful for the hospitality and collaboration. She also acknowledges the ARC Centre of Excellence for Mathematical and Statistical Frontiers for partially supporting her research visit at this University.



\addcontentsline{toc}{section}{References}
\bibliographystyle{plain}


\end{document}